\documentstyle{amltd2004}
\begin{document}
\annalsline{155}{2002}
\received{February 13, 2000}
\startingpage{405}
\def\bye{\end{document}}
 \font\tenrm=cmr10
\def\ritem#1{\item[{\rm #1}]}
\input boxedeps.tex 
\SetepsfEPSFSpecial 
\HideDisplacementBoxes
\def\figin#1#2{$$
 {\BoxedEPSF{xy#1.eps scaled
#2}%
}%
$$
\noindent}
\input amssym.def
\input amssym.tex

\def\joinrel{\mathrel{\mkern-4mu}}
\def\relbar{\mathrel{\smash-}}
\def\lrar{\relbar\joinrel\relbar\joinrel\rightarrow}

\catcode`\@=11
\font\twelvemsb=msbm10 scaled 1100
\font\tenmsb=msbm10
\font\ninemsb=msbm10 scaled 800
\newfam\msbfam
\textfont\msbfam=\twelvemsb  \scriptfont\msbfam=\ninemsb
  \scriptscriptfont\msbfam=\ninemsb
\def\msb@{\hexnumber@\msbfam}
\def\Bbb{\relax\ifmmode\let\next\Bbb@\else
 \def\next{\errmessage{Use \string\Bbb\space only in math
mode}}\fi\next}
\def\Bbb@#1{{\Bbb@@{#1}}}
\def\Bbb@@#1{\fam\msbfam#1}
\catcode`\@=12

 \catcode`\@=11
\font\twelveeuf=eufm10 scaled 1100
\font\teneuf=eufm10
\font\nineeuf=eufm7 scaled 1100
\newfam\euffam
\textfont\euffam=\twelveeuf  \scriptfont\euffam=\teneuf
  \scriptscriptfont\euffam=\nineeuf
\def\euf@{\hexnumber@\euffam}
\def\frak{\relax\ifmmode\let\next\frak@\else
 \def\next{\errmessage{Use \string\frak\space only in math
mode}}\fi\next}
\def\frak@#1{{\frak@@{#1}}}
\def\frak@@#1{\fam\euffam#1}
\catcode`\@=12
\newfont{\bm}{msbm10}

\newcommand{\semi}{\mbox{{\bm \symbol{111}}}} 

 
\newcommand{\ncov}{\langle n \rangle}
\newcommand{\onecov}{\langle 1 \rangle}
\newcommand{\map}{\textstyle\mathop{\rm map}}
\newcommand{\pho}{\textstyle\mathop{\rm Ho}_*}
\newcommand{\qui}{\textstyle\mathop{\rm \bf A}}
\newcommand{\Fib}{\textstyle\mathop{\rm Fib}}
\newcommand{\pG}{\textstyle\mathop{\rm polyGEMs}}
\newcommand{\opG}{\textstyle\mathop{\rm polyGEMs^{ori}}}
\newcommand{\pGf}{\textstyle\mathop{\rm polyGEMs^{ft}}}
\newcommand{\pkt}{*}
\newcommand{\s}{\Sigma}
\newcommand{\ho}[1]{\textstyle\mathop{\rm Ho}(#1)}
\newcommand{\Vect}{\textstyle\mathop{\rm Vect}}
\newcommand{\sd}{\textstyle\mathop{\rm sd}}
\newcommand{\di}{\textstyle{\rm d}}  
\newcommand{\dio}{\textstyle\mathop{\rm \bar{d}}}
\newcommand{\K}{\textstyle\mathop{\rm \bf K}}
\newcommand{\sta}{\textstyle\mathop{\rm \bf star}}
\newcommand{\link}{ \textstyle\mathop{\rm \bf link}}
\newcommand{\Top}{ \textstyle\mathop{\rm Top}}
\newcommand{\nstar}{\textstyle\mathop{\rm star}}
\newcommand{\pequiv}{\stackrel{\F_p}{\cong}}
\newcommand{\Cone}{\textstyle\mathop{\rm Cone}}

\newcommand{\Syl}{\textstyle\mathop{\rm Syl}}
\newcommand{\ch}{\textstyle\mathop{\rm char}}
\newcommand{\colim}{\textstyle\mathop{\rm colim}}
\newcommand{\aut}{\textstyle\mathop{\rm aut}}
\newcommand{\Aut}{\textstyle\mathop{\rm Aut}}
\newcommand{\Rep}{\textstyle\mathop{\rm Rep}}
\newcommand{\Out}{\textstyle\mathop{\rm Out}}
\newcommand{\End}{\textstyle\mathop{\rm End}}
\newcommand{\rk}{\textstyle\mathop{\rm rk}}
\newcommand{\Spec}{\textstyle\mathop{\rm Spec}}
\newcommand{\Ext}{\textstyle\mathop{\rm Ext}}
\newcommand{\Hom}{\textstyle\mathop{\rm Hom}}
\newcommand{\mor}{{\rm Mor}}
\newcommand{\Tor}{\textstyle\mathop{\rm Tor}}
\newcommand{\id}{\textstyle\mathop{\rm id}}
\newcommand{\res}{\textstyle\mathop{\rm res}}
\newcommand{\Der}{\textstyle\mathop{\rm Der}}
\newcommand{\diag}{\textstyle\mathop{\rm diag}}
\newcommand{\Rmod}{R\mbox{-}\textstyle\mathop{\rm mod}}
\newcommand{\Zpmod}{\Z_{(p)}\mbox{-}\textstyle\mathop{\rm mod}}
\newcommand{\Zlpmod}{\Z_{(p)}\mbox{-}\textstyle\mathop{\rm mod}}
\newcommand{\Zmod}{\Z\mbox{-}\textstyle\mathop{\rm mod}}
\newcommand{\Lmod}{\Lambda\mbox{-}\textstyle\mathop{\rm mod}}
\newcommand{\Zpm}[1]{{{\Z_{(p)}}{#1}\mbox{-}\textstyle\mathop{\rm mod}}}
\newcommand{\Rm}[1]{{{R}{#1}\mbox{-}\textstyle\mathop{\rm mod}}}
\newcommand{\Lm}[1]{{{\Lambda}{#1}\mbox{-}\textstyle\mathop{\rm mod}}}
\newcommand{\Op}{{\bf O}_p(G)}
\newcommand{\pR}{{\cal B}_p(G)}
\newcommand{\pRe}{{\cal B}_p^e(G)}
\newcommand{\pS}{{\cal S}_p(G)}
\newcommand{\pSe}{{\cal S}_p^e(G)}
\newcommand{\pA}{{\cal A}_p(G)}
\newcommand{\pD}{{\cal D}_p(G)}
\newcommand{\pDeM}{{\cal D}_{p,M}^e(G)}
\newcommand{\pDM}{{\cal D}_{p,M}(G)}
\newcommand{\pI}{{\cal I}_p(G)}
\newcommand{\pH}{{\cal H}_p(G)}

\newcommand{\holim}{\textstyle\mathop{\rm holim}}

\newcommand{\fF}{{\frak F}}
\newcommand{\cF}{{\cal F}}

\newcommand{\cC}{{\cal{C}}}
\newcommand{\hi}{{\rm ht}}
\newcommand{\bO}{{\bf O}}
\newcommand{\coker}{{\rm coker}}
\newcommand{\op}{{\rm op}}

\newcommand{\Cat}{{\rm {\bf Cat}}}
\newcommand{\Spaces}{{\rm {\bf Spaces}}}

\newcommand{\St}{{\rm St}}
\newcommand{\SSt}{{\rm {\Bbb S}t}}
\newcommand{\pfrep}{{\rm fusion-Rep}}

\newcommand{\g}{{\frak g}}
\newcommand{\h}{{\frak h}}
\newcommand{\Na}{{\frak a}}
\newcommand{\r}{{\frak r}}
\newcommand{\n}{{\frak n}}
\newcommand{\Nb}{{\frak b}}
\newcommand{\Ns}{{\frak s}}
\newcommand{\Nht}{\textstyle\mathop{\rm ht}}
\newcommand{\fraks}{{\frak s}}
\newcommand{\gl}{{\frak gl}}
\newcommand{\so}{{\frak so}}
\newcommand{\Nsl}{{\frak sl}}
\newcommand{\Nsp}{{\frak sp}}
\newcommand{\GL}{{\rm GL}}
\newcommand{\SP}{{\rm SP}}
\newcommand{\SO}{{\rm SO}}
\newcommand{\SL}{{\rm SL}}
\newcommand{\Lie}{{\rm Lie}}
\newcommand{\Tr}{{\rm Tr}}
\newcommand{\Mat}{{\rm Mat}}
\newcommand{\ad}{{\rm ad}}

\newcommand{\luft}{\medskip\par\noindent}
\newcommand{\QED}{\hfill {\bf $\square$}\luft}
\newcommand{\start}[1]{\begin{#1}\em}
\newcommand{\slut}[1]{\end{#1}}

\newcommand{\C}{{\bf {C}}}
\newcommand{\Q}{{\bf {Q}}}
\newcommand{\R}{{\bf {R}}}
\newcommand{\N}{{\bf {N}}}
\newcommand{\F}{{\bf {F}}}
\newcommand{\Z}{{\bf {Z}}}
\newcommand{\E}{{\bf {E}}}
\newcommand{\k}{{\bf{k}}}

\newcommand{\T}{{\bf {T}}}
\newcommand{\D}{{\bf D}}
\newcommand{\CP}{{\bf{CP}}}
\newcommand{\RP}{{\bf{RP}}}
\newcommand{\QP}{{\bf{QP}}}
\newcommand{\p}{{\bf{P}}}

\newcommand{\cK}{{\cal{K}}}
\newcommand{\nil}{{\cal{N}}{\it il} }
\newcommand{\bnil}{\overline{{\cal{N}}{\it il}} }
\newcommand{\B}{{\cal{B}} }
\newcommand{\cU}{{\cal{U}}}
\newcommand{\cG}{{\cal{G}}}
\newcommand{\PS}{{\cal{PS}}}
\newcommand{\cE}{{\cal{E}}}
\newcommand{\A}{{\bf A}}
\newcommand{\Sq}{\textstyle\mathop{\rm Sq}}

\newcommand{\ddx}{\frac{d}{dx}}
\newcommand{\ddxk}{\frac{d^k}{dx^k}}
\newcommand{\lsim}{\stackrel<\sim}
\newcommand{\gsim}{\stackrel>\sim}
\newcommand{\into}{\hookrightarrow}
\newcommand{\onto}{\twoheadrightarrow}


\newcommand{\func}[3]{\mbox{$#1\colon #2\rightarrow #3$}}
\newcommand{\beq}{\begin{eqnarray*}}
\newcommand{\eeq}{\end{eqnarray*}}
\newcommand{\pil}[1]{\stackrel{#1}{\rightarrow}}
\newcommand{\bil}[1]{\stackrel{#1}{\leftarrow}}
\newcommand{\bto}{\leftarrow}
\newcommand{\ul}{\underline}
\newcommand{\tuborg}{\left\{\begin{array}{ll}}
\newcommand{\sluttuborg}{\end{array}\right.}


\title{Higher limits via subgroup complexes} 

 \author{Jesper Grodal}
 \institutions{ Massachusetts Institute of Technology,
        Cambridge, MA\\
{\eightpoint {\it  Current address\/}:} University of Chicago, Chicago, IL\\
{\eightpoint {\it E-mail address\/}: jg@math.uchicago.edu}}

\section{Introduction}

Abelian group-valued functors on the orbit category of a finite group $G$
are ubiquitous both in group theory and in much of homotopy theory.
In this paper we give a new finite model for $\lim^*_\D F$, the higher
derived functors of the inverse limit functor, for a general functor $F :
\D^{\op} \to \Zpmod$, where $\D$ is the $p$-orbit category of a finite group
$G$ (consisting of transitive $G$-sets
 with isotropy groups, $p$-groups and equivariant maps between them), or one
of several other commonly studied
categories associated to $G$.

In homotopy theory such functors arise in connection with homology
decompositions of the classifying space $BG$ of $G$, where $\D$ serves as
indexing category for the decomposition. The groups $\lim^*_\D F$ then occur
in the $E_2$-term of several fundamental spectral sequences such as the
Bousfield-Kan homotopy (or cohomology) spectral sequence of a homotopy limit
(or colimit)
$$ E_2^{i,j} = \lim_\D{}^{-i} \pi_j(X(d)) \Rightarrow \pi_{i+j}(\holim_\D
X(d)),$$
where $X$ is a functor $\D^{\op} \to  \Spaces_{*}$ (\cite{BK72} or unpointed
spaces \cite{bousfield89}). In
applications, $\holim_\D X(d)$ will be a space of significant topological
interest such as the space of maps from
$BG$ to some other given space. It is hence of great importance to be able
to calculate these higher limits for
calculating maps between classifying spaces  \cite{JMO92}, \cite{JO96},
\cite{vectmanus}, finding homology
decompositions
\cite{JM92}, \cite{dwyer97}, \cite{dwyer98}, and settling uniqueness issues
\cite{notbohm94}.

In group theory, $F$ is typically a Mackey functor such as a group
cohomology. Here one can often show that the higher
limits vanish for $*>0$ for general reasons (cf.\ \cite{JM92}, \cite{webb99}
and Section~\ref{exactseqsec}). In this
seemingly trivial case the general formulas for $\lim^*_\D F$ in fact
specialize to give interesting expressions, in
the form of finite exact sequences, for the `group of stable elements'
$\lim^0_\D F$, which is equal, or closely
related, to $F(G)$.
\vfill
\footnoterule
\vglue4pt
{\ninerm   2000 {\nineit Mathematics Subject Classification}. Primary:
55R35; Secondary: 18G10, 55R37, 20J06,
55N91.} 
\pagebreak

By a collection we mean a set $\cC$ of subgroups in $G$, closed under
conjugation in $G$. We can view a collection as a poset under inclusion of
subgroups, and the nerve $|\cC|$ of this poset is what is called the {\it
subgroup complex} associated to the collection $\cC$.  Obviously it is a
{\it finite} (ordered) $G$-simplicial complex. Subgroup complexes play an
important role in finite group theory, where they are thought of as
generalized buildings (or geometries) of $G$.

Let $\bO_\cC$ denote the $\cC$-orbit category, that is, the full subcategory
of the orbit category $\bO(G)$ with
objects the transitive $G$-sets with isotropy groups in $\cC$. This category
arises in connection with the `subgroup
homology decomposition'
\cite{JMO92}, \cite{dwyer97}, \cite{dwyer98}. (We will also consider the
(opposite) $\cC$-conjugacy category $\A_\cC^{\op}$, associated with the
`centralizer decomposition', and the orbit simplex category $(\sd \cC)/G$
associated with the `normalizer decomposition'.)

We show that for many collections $\cC$, such as the collection of all
nontrivial $p$-subgroups, higher limits $\lim_\D^*F$ over any of the three
types of categories can be interpreted as Bredon equivariant cohomology of
$|\cC|$ with values in a $G$-local coefficient system
(or presheaf) constructed from $F$. Here a $G$-local coefficient system is
understood in the sense of Ronan-Smith
\cite{RS85}, \cite{RS86}, \cite{benson98}, which is close to the original
definition of Bredon \cite{bredon67}, but is more general
than what is usually called a (Bredon) coefficient system, and this is a key
generalization which the result depends
upon. This result gives a finite model for the higher limits, breaking the
complexity of the orbit category up into two
parts, separating the poset structure from the $G$-conjugation, which can in
some sense be viewed as a topological
analog of Alperin's fusion theorem \cite{alperin67}. Computationally the
model is a large improvement over earlier
Bredon cohomology models for the higher limits over $\bO_\cC^{\op}$ and
$\A_\cC$ which have been infinite
\cite{mislin89}, \cite{JM89}, \cite{JMO92}, \cite{webb99}. It furthermore
reveals a close connection between higher limits and modern
finite group theory, and gives an easy systematic way to see which subgroups
can be omitted from $\cC$ without
changing the higher limits.

We now give a more detailed description of our results in the case of
$\bO_\cC^{\op}$. (We obtain similar results for $\A_\cC$ and $((\sd
\cC)/G)^{\op}$.)
In order to do so, we first need to introduce some notation.
Let $\pS$ and $\pR$ denote the collection of {\it nontrivial} $p$-subgroups
and $p$-radical subgroups respectively.
 Recall\break that
a {\it $p$-radical subgroup} of $G$ is a $p$-subgroup $P$ with the property
that\break $O_p(NP/P)=~e$. (Here,
$O_p(\cdot)$ means largest normal $p$-subgroup.) We use the notation ${\cal
S}^e_p(G)$ etc.\ to denote the same
posets as above, but now not excluding the trivial subgroup $e$.
\vglue2pt

Let $\dio \cC$ denote the opposite category of the simplex category of
$|\cC|$. The objects are the simplices of $|\cC|$, i.e., chains of subgroups
$H_0 \leq \cdots \leq H_n$, with $H_i \in \cC$, and the morphisms are
generated by the opposite of face and degeneracy maps. Let $(\dio \cC)_G$ be
the category that has the same objects as $\dio \cC$ but where a morphism
$\sigma \to \sigma'$ is a pair $(f: g\sigma \to \sigma',g)$, where $f$ is a
morphism in  $\dio \cC$ and $g \in G$, and morphisms
are composed via $(f,g) \circ (f',g') = (f \circ (gf'),gg')$ (see
Section~2.1).
By a $G$-local coefficient system  (over $\Z_{(p)}$) on $|\cC|$ we mean a
(covariant) functor $\cF: (\dio \cC)_G \to \Zpmod$. For a functor on
$\bO(G)$ and its subcategories we write $F(H)$ for $F$ evaluated on $G/H$.

Given a functor $F : \bO_\cC^{\op} \to \Zpmod$ we equip $|\cC|$ with the
$G$-local coefficient system $\cF$  given by precomposing
it  with the canonical map $(\dio \cC)_G \to (\cC^{\op})_G \to
\bO_\cC^{\op}$ (see Remark \ref{grothrem}). On objects,
we have $\cF(H_0 \leq \cdots \leq H_n) = F(H_0)$.

\proclaim{Theorem} \label{thm1} Let $\cC$ be any collection of $p$\/{\rm
-}\/subgroups of a finite group $G${\rm ,}
 closed under passage to $p$\/{\rm -}\/radical overgroups{\rm .}
 Let  $F: \bO_\cC^{\op} \to \Zpmod$ be any functor{\rm ,} and let $\cF$ be
the $G$\/{\rm -}\/local coefficient system on $|\cC|$ given by precomposing
$F$ with
 $(\dio \cC)_G \to (\cC^{\op})_G \to
\bO_\cC^{\op}${\rm .} Then
$$\lim_{\bO_\cC}{}^*F \cong H^*_G(|\cC|;\cF).$$

Furthermore{\rm ,} for any collection $\cC'$ satisfying $\cC \cap \pRe
\subseteq \cC' \subseteq \cC${\rm ,} the inclusion
$|\cC'| \to |\cC|$ induces an isomorphism $ H^*_G(|\cC|;\cF) \pil{\simeq}
H^*_G(|\cC'|;\cF)${\rm ,} where $\cF$ is a
coefficient system on $|\cC'|$ by restriction{\rm .}

If $F$ is a functor concentrated on conjugates of a single subgroup $P${\rm
,}
 then$$ \lim_{\bO_\cC}{}^i F \cong H^{i-1}( \Hom_{NP/P}(\St_{*}(NP/P),F(P)))
.$$
\endproclaim 

Here $H^*_G(|\cC|;\cF)$ denotes Bredon cohomology $H^*_G(|\cC|;\cF) =
H(C^*(|\cC|;\cF)^G)$ (see Definition \ref{bredoncohom}), which of course is
zero above the dimension of $|\cC|$, and $\St_*(W)$ denotes the Steinberg
complex of $W$, the reduced normalized chain complex $\tilde C_*(|{\cal
S}_p(W)|;\Z_{(p)})$.
The result that one only needs the $p$-radical subgroups for calculating the
higher limits is due to Jackowski-McClure-Oliver \cite{JMO92}.
They used a Bredon cohomology model (\cite{mislin89}, \cite{JM89}, and
Prop.\ \ref{prop2}) but the space was
infinite, and hence hard to use directly for calculations.

By the work of Webb (\cite{webb87arcata}, \cite{webb91}, and see
Theorem~\ref{webbthm}) it is known that the
Steinberg complex over the $p$-adic integers $\Z_p$ is $G$-chain homotopic
to a (unique minimal) direct summand
which is a chain complex of projective\break $\Z_pG$-modules; we denote this
by $\SSt_*(G)$. This chain complex has
been recognized by group theorists as a worthy generalization of the
Steinberg module to arbitrary finite groups
\cite{alperin90}, \cite{webb87arcata}, \cite{webb91} and we study its
properties in Section~\ref{st}. (If $G$ is a
finite group of Lie type of characteristic $p$ then $\SSt_*(G)$ is just the
ordinary Steinberg module concentrated in
degree one less than the semisimple rank of $G$.)  Theorem~\ref{thm1}
combined with (mostly straightforward)
properties of the Steinberg complex give significant generalizations of the
properties of higher limits over functors
concentrated on one conjugacy class, so called {\it atomic} functors, found
in \cite{JMO92}.

In homotopy theory, the functors which occur in practice often have the
property that the centralizer $CP$ of $P$ acts trivially on $F(P)$ (i.e.,
the action of $NP/P$ on $F(P)$ factors through $NP/PC(P)$). In this case one
can remove additional subgroups from the collection. Let $\pD$ be the
collection consisting of the subgroups which are {\it $p$-centric} (i.e.,
the center $ZP$ of $P$ is a Sylow $p$-subgroup in $CP$), and such that
$O_p(NP/PC(P)) = e$. These subgroups are necessarily $p$-radical subgroups
and, like the $p$-radical groups, play a role in Alperin's theory of weights
for arbitrary finite groups, where they are the $p$-radical subgroups for
the principal block  \cite[\S 3]{alperin87}. For this reason we call the
subgroups in $\pD$ the {\it principal} $p$-radical subgroups.

\proclaim{Theorem} \label{pcentriccor} \hskip-8pt
Let $\cC$ be any collection of $p$\/{\rm -}\/subgroups of a finite
group~$G${\rm ,}
 closed under passage to $p$\/{\rm -}\/radical overgroups{\rm .}
Suppose that $F$ is a functor $\bO_\cC^{\op} \to \Zpmod$ such that for all
$P\in \cC${\rm ,}
 $CP$ acts trivially on $F(P)${\rm .}
Then{\rm ,} for any collection $\cC'$ satisfying $\cC \cap \pD \subseteq
\cC' \subseteq \cC${\rm ,}
$$\lim_{\bO_\cC}{}^*F \cong \lim_{\bO_{\cC'}}{}^*F \cong H^*_G(|\cC'|;\cF) ,
$$
where $\cF$ is the $G$\/{\rm -}\/local
coefficient system on $|\cC'|$ induced via  $(\dio \cC')_G \to (\dio \cC)_G
\to
(\cC^{\op})_G \to \bO_\cC^{\op}${\rm .}

If $F$ is a functor concentrated on conjugates of a single subgroup $P \in
\cC \cap \pD${\rm ,}
 then$$ \lim_{\bO_\cC}{}^i F \cong H^{i-1}(
\Hom_{NP/PC(P)}(\St_{*}(NP/PC(P)),F(P))) .$$
\endproclaim 

The formulas for higher limits of functors concentrated on one conjugacy
class easily show that the questions of minimality of the collections $\cC
\cap \pRe$ in Theorem \ref{thm1} and $\cC \cap \pD$ in Theorem
\ref{pcentriccor} both, in a certain sense, are equivalent to a version of
Quillen's conjecture on the contractibility of $|{\cal S}_p(W)|$ (see
Remarks \ref{quillenrem1} and \ref{quillenrem2}).

The more explicit results in Theorems~\ref{thm1} and~\ref{pcentriccor} for
atomic functors provide a computation of
the $E_1$-term of a spectral sequence converging to the higher limits of a
general functor $F$ obtained by filtering
$F$ in such a way that the associated quotients are atomic functors. By a
$G$-invariant height function $\hi$ on a
collection $\cC$ we mean a strictly increasing map of posets $\cC/G \to \Z$.

\proclaim{Theorem}\label{sscor} \hskip-8pt Let $\cC$ be any collection of
$p$\/{\rm -}\/subgroups of a
finite group~$G${\rm ,} closed under passage to $p$\/{\rm -}\/radical
overgroups{\rm .}
For any functor $F: \bO_\cC^{\op} \to \Zpmod$ there is  a spectral sequence
\begin{eqnarray*}
E_1^{i,j} &=& \bigoplus_{[P] \in \cC'/G,
\hi(P)=-i}H^{i+j-1}(\Hom_{NP/P}(\St_*(NP/P),F(P)))
\\[5pt]
&\Rightarrow& \lim_{\bO_\cC}{}^{i+j}F
\end{eqnarray*}
where $\cC' = \cC \cap \pRe$ and $\hi$ is some chosen $G$\/{\rm
-}\/invariant height function on~$\cC'${\rm .}
The differentials are cohomological differentials $d_n: E_n^{i,j} \to
E_n^{i+n,j-(n-1)}${\rm .}

If for all $P \in \cC${\rm ,} $CP$ acts trivially on $F(P)${\rm ,}
 then  $\cC'$ may be replaced by $ \cC \cap \pD$ and $NP/P$ by
$NP/PC(P)${\rm .}
\endproclaim 

We can always normalize a height function such that for a Sylow $p$-subgroup
$P$, $\hi(P) = 0$, and
then this spectral sequence occupies a triangle in the `eighth octant'.
Often the most convenient $G$-invariant height function on a collection
$\cC$ is the `minimal' one given by $\hi(Q)~=~-\dim |\cC_{\geq Q}|$.
If $G$ is a finite group of Lie type of characteristic $p$ the poset  $\pRe$
has a relatively simple structure (it is Cohen-Macaulay). In this case, if
we choose the right height function, the $E_1$-term is concentrated on the
horizontal axis, and hence yields a chain complex, with an explicitly
described differential, for calculating the higher limits (see
Theorem~\ref{liecor}).

\demo{Applications to group cohomology and other Mackey functors}
In cases where $\lim^i$ vanishes for $i>0$, the chain complex defining
Bredon cohomology in Theorem~\ref{thm1}
 produces finite exact sequences starting with $\lim^0_\D F$. (We get
similar exact sequences when treating higher
limits over $\cC$-conjugacy categories and orbit simplex categories.) For
group cohomology this gives exact
sequences expressing $H^n(G;M)$ in terms of $H^n(H;M)$ for proper subgroups
$H$ of $G$.
\enddemo

\proclaim{Theorem}  Let $G$ be a finite group{\rm ,} and let $H^n(-)$ denote
group cohomology with coefficients in
$\Z_{(p)}${\rm .} There are  the following three exact sequences coming from
acyclic complexes
\begin{eqnarray*}
 0 \to H^n(G) 
&\hskip-7pt\to &\hskip-7pt\oplus_{[P] \in \pD/G} H^n(P)^{NP}\\
&\hskip-7pt\to &\hskip-7pt\oplus_{[P_0 < P_1] \in |\pD|_1/G}
H^n(P_0)^{NP_0\cap NP_1} \to \cdots
\to 0, \\
  0 \to H^n(G)&\hskip-7pt\to &\hskip-7pt\oplus_{[V] \in \pA/G} H^n(CV)^{NV}
\\
&\hskip-7pt\to &\hskip-7pt\oplus_{[V_0 < V_1] \in |\pA|_1/G}
H^n(CV_1)^{NV_0\cap NV_1} \to \cdots \to 0,\\
 0 \to H^n(G) &\hskip-7pt\to &\hskip-7pt\oplus_{[P] \in \pD/G} H^n(NP)\\
&\hskip-7pt\to &\hskip-7pt\oplus_{[P_0<P_1]
\in |\pD|_1/G}  H^n(NP_0\cap NP_1) \to \cdots \to 0.
\end{eqnarray*}
\endproclaim 
Here $\pA$ denotes the collection of nontrivial elementary abelian
$p$-subgroups, that is subgroups of the form $(\Z/p)^k$, for some $k>0$, and 
$|\cC|_k/G$ denotes $G$-conjugacy classes of nondegenerate $k$-simplices in 
$|\cC|$.
The last exact sequence is a version of the celebrated exact sequence of
Webb, but obtained for the small collection $\pD$, whereas the rest appear
to be new.
We explain these sequences as part of a general framework in
Section~\ref{exactseqsec}, and also give related exact
sequences.\pagebreak

The spectral sequence in Theorem \ref{sscor} can also be used to give a
useful description of $\lim^0_{\bO_\cC} F$, regardless of the vanishing of
the higher limits, simply by looking at what happens at the $(0,0)$-spot. In
this way we recover expressions for the stable elements, usually proved
using various versions of Alperin's fusion theorem, via homotopy theoretic
considerations (see  Section~\ref{stableelements}).

\demo{Applications to homology decompositions}
A homology decomposition is said to be {\it sharp} \cite{dwyer98} if the
Bousfield-Kan cohomology spectral sequence of the homotopy colimit collapses
at the $E_2$-term onto the vertical axis. (This can also be formulated in
other equivalent ways---see Section~\ref{sharpsec}.)
Our general formulas for higher limits can be used to obtain new sharpness
results for homology decompositions. It is easy to see that the collection
$\pSe$ is sharp with respect to both the subgroup, centralizer and
normalizer homology decompositions.
Theorem~\ref{pcentriccor}, for $F$ equal to mod $p$ group cohomology, shows
that the $E_2$-term of the cohomology spectral sequence associated to the
subgroup decomposition using all $p$-subgroups remains unchanged when
all subgroups not in $\pD$ are removed. This, together with the analogous
result for the normalizer decomposition
(Theorem~\ref{normsharp}), yields
\enddemo

\proclaim{Theorem} \label{introsharpnessresult} The collection $\pD$ of
principal $p$\/{\rm -}\/radical
subgroups gives a subgroup and normalizer sharp {\rm mod} $p$ homology
decomposition of $BG$ {\rm (}\/and is in
 particular
ample\/{\rm ).}
\endproclaim 

(See Theorem~\ref{sharpnessresult} for a more elaborate version of this
result.)
The subgroup sharpness result improves a result of Dwyer \cite{dwyer97},
\cite{dwyer98} that the collection of
$p$-radical and $p$-centric subgroups is subgroup sharp. The normalizer
sharpness result improves a result of Webb
\cite{webb87arcata} that the collection of $p$-radical subgroups is
normalizer sharp and verifies a suspicion of
Smith-Yoshiara \cite{SY97} that the collection of $p$-radical and
$p$-centric subgroups should be normalizer sharp.

\demo{Related work}
This work grew out of a need to calculate higher limits\break (initially
using \cite{JMO92}) in order to compute
monoids of vector bundles over classifying spaces completed at a prime $p$,
as a part of a larger program to
understand the `homotopy representation theory' of finite, compact Lie, or
even $p$-compact groups
\cite{dwyer98icm}, at the prime
$p$. This motivating application of our results will be published separately
\cite{vectmanus}.

The relationships between groups, posets, and spaces which underlie
everything in this paper date  back to the
fundamental work of Brown \cite{brown75} and Quillen \cite{quillen78}. We
assume some familiarity with the
methods of \cite{quillen78}. The approach taken in this paper can be seen as
a continuation along the lines of the
papers of Dwyer \cite{dwyer97}, \cite{dwyer98}, which we highly recommend to
the reader, although this paper can be
read independently of those.

The methods of this paper are also closely related to
a `Lie theoretic' approach to local representation theory of finite groups
advocated by Alperin
\cite{alperin87}, \cite{alperin90} and Webb \cite{webb87arcata} amongst
others, and likewise have ties to early
ideas of Puig \cite{puig76} who examines the basic but yet elusive question
of what is a $p$-local group?

Results relating higher limits over various categories to higher limits over
orbit simplex categories have been found independently by J. S\l omi\'nska
\cite{slominska96}.  More precisely, using for instance Proposition
\ref{normprop} (to remove a subdivision), one sees that 
 \cite[Prop.\ 2.9]{slominska96} (which is formulated in a more abstract
setup) corresponds to Corollary
\ref{cor1'} which is an important step in the proof of Theorem \ref{thm1}. I
am grateful to S. Jackowski for informing
me of her unpublished work and for providing me with a copy. After the
submission of the present paper Jackowski and
S\l omi\'nska have finished the paper \cite{JS01} which includes S\l
omi\'nska's unpublished results.
\enddemo

{\it Organization of the paper}.
In Section~\ref{prelim} we prove various facts about higher limits and
Bredon cohomology needed in the later sections. In Section~\ref{orbitsec} we
prove our general results on higher limits over orbit categories,
several times referring to results from Section~\ref{st}.
Section~\ref{liesec} elaborates on these results in the
special case where $G$ is a finite group of Lie type of characteristic $p$. 
In Section~\ref{st} we
examine the structure of the Steinberg
complex
$\St_*(G)$ for an arbitrary finite group $G$. Sections~\ref{conjugacysec}
and \ref{simplexsec} treat higher limits
over conjugacy categories and orbit simplex categories respectively.  In the
last three sections we look at special
cases of the general results, focusing on the case where the higher limits
in fact vanish. In Section~\ref{exactseqsec}
we explain how one can obtain exact sequences in the case where all higher
limits vanish, and provide some examples
when this is satisfied.  In Section~\ref{sharpsec} we explain how the
general results lead to new sharpness results
on homology decompositions. Finally, in Section~\ref{stableelements} we
explain how our results relate to the
calculation of stable elements and Alperin's fusion theorem.

\vfill
{\it Notation and conventions}.
Throughout, $G$ denotes a finite group and $p$ a fixed prime.
This paper is written simplicially so  that the word space, without any
prefix, means a simplicial set. Hence the nerve
construction $|\cdot|$ is a functor from categories to simplicial sets. We
remark in this connection  that for\break
$G$-categories, taking $G$-fixpoints commutes with taking nerves. We let $R$
denote an arbitrary commutative
unital ring. In most of this paper we however work over $\Z_{(p)}$, and we
adopt the convention that unless otherwise
indicated, $\Z_{(p)}$ is taken as the ground ring. We denote the normalizer
and centralizer of $H$ in $G$ by $NH$ and
$CH$ respectively, and write $WH$ for the Weyl group $NH/H$. The word
`functor' without any prefix means a
covariant functor. The suspension $\Sigma X$ of a $G$-space $X$ should be
understood as ${\Cone}X/X$, where
${\Cone}X$ here is the {\it unreduced} cone on $X$. This is a pointed
$G$-space, with $G$-basepoint $[X]$.
\eject

\demo{Acknowledgments}
I would like to thank Bill Dwyer, Dan Kan, Bob Oliver, and Peter Symonds for
helpful conversations. Exchange of
calculations with Bob Oliver led to improvements in the statement of the
results in Section~\ref{st}. I thank Kasper
Andersen and Peter Webb for making useful comments on a preliminary version
of this manuscript, and Lars
Hesselholt for providing helpful advice on the exposition. Finally, I would
like to thank my advisor Haynes Miller for
many helpful conversations and comments and much encouragement.

\section{Preliminaries} \label{prelim}

In this section we explain various facts about equivariant cohomology and
higher limits. All of the results explained in this section are
 well-known. 

\demo{{\rm 2.1.} Equivariant cohomology and local coefficient systems}
\label{equiv}
Let $X$ be a space. Define the simplex category (or division \cite{DK83})
$\di X$ of $X$ to be the category with objects the simplices of $X$ and
morphisms generated by face and degeneracy maps. In other words the objects
are maps $\sigma: \Delta[n] \to X$, where $\Delta[n]$ is the standard
$n$-simplex, and a morphism $\sigma \to \sigma'$ is given by a commutative
triangle
\figin{1}{1000}

\noindent 
Often we need the opposite category of $\di X$, which we denote $\dio X$.
If $X$ is the nerve of a small category $\D$, then we abbreviate $\di|\D|$
to $\di\D$. 
Define a {\it local coefficient system} on $X$ (over $R$) to be a functor
$\cF:\dio X \to \Rmod$.

For any local coefficient system $\cF$ we can define the simplicial cochain
complex $(C^*(X;\cF),\delta)$ of $X$ with coefficients in $\cF$ by
 setting $$C^n(X;\cF) = \prod_{\sigma \in X_n}\cF(\sigma) \hbox{ and }
 (\delta f)(\sigma) = \sum_i (-1)^i\cF(d_i^{\op})(f(d_i \sigma)).$$
We define $H^*(X;\cF)$ to be the homology of this
cochain complex (cf.\ also \cite[p.\ 42ff]{godement58}). As usual we may as
well use the {\it normalized} simplicial
cochain complex for calculating $H^*(X;\cF)$ (cf.\ \cite[\S 8.3]{weibel94}).

Note that we have a natural transformation $\dio \D \to \D$ on objects given
by $(x_0 \to \cdots \to x_n) \mapsto x_n$.
We also have a natural transformation $\dio\D \to \D^{\op}$ on objects given
by $(x_0 \to \cdots \to x_n) \mapsto x_0$. If $F: \D \to \Rmod$ is a
covariant (resp.\ contravariant) functor then $F$ induces a local
coefficient system $\cF$ on $|\D|$ via  $\dio \D \to \D$ (resp.\ $\dio \D
\to \D^{\op}$).
In both the co- and contravariant case we denote the cohomology groups by
$H^*(\D;F)$, i.e., we omit indicating taking nerves, and indicate that the
local coefficient system comes from a functor on $\D$ by keeping the
notation $F$ instead of $\cF$.

We now want to define Bredon equivariant cohomology. For our purposes it is
most natural to use a notion of a $G$-local coefficient system which is more
general than the one given in e.g.\  \cite{bredon67}, but which agrees with
the usage in group theory (cf.\ \cite[Ch.\ 7]{benson98}). Roughly speaking,
we want to allow different elements of $G$ to induce different automorphisms
of $C^*(X;\cF)$ even though they act the same way on $X$.
For this purpose we introduce the following construction.
Let $\D$ be a small $G$-category, and let $\alpha$ be the functor $\alpha: G
\to \Cat$ which defines this
$G$-category. ($G$ is viewed as a category with one object, and $\Cat$
denotes the category of small categories.)
Note that if $\D$ is a $G$-category, then so is $\D^{\op}$, via $\op: \Cat
\to \Cat$. Set $\D_G$ equal to the
Grothendieck construction of $\alpha$ \cite{thomason79}, i.e., the category
which has the same objects as $\D$ but
where the morphisms $x \to y$ in $\D_G$ are given by pairs $(f: gx \to
y,g)$, where $f$ is a morphism in $\D$ and $g
\in G$, which compose via $(f',g') \circ (f,g) = (f' \circ (g'f),g'g)$. Note
that this can be thought of as a sort of
semi-direct product of $\D$ with $G$, and just amounts to enlarging the
morphisms in $\D$ with the elements of $G$
in the freest possible way.

\advance\theoremcount by 1

\numbereddemo{{R}emark} \label{grothrem}
Let $\D = \cC$ for some collection $\cC$. The categories $\cC_G$,
$\bO_\cC^{\op}$ and $\A_\cC$ can all be viewed as having the same objects,
namely the objects of $\cC$. The morphisms from $H$ to $H'$ in $\cC_G$ can
be identified with the elements of $g \in G$ such that $gHg^{-1}\leq H'$.
The morphisms from  $H'$ to $H$ in $\bO_\cC^{\op}$ may be identified with
the elements $[g] \in G/H'$ such that $gH'{g^{-1}} \geq H$ and the morphisms
from $H$ to $H'$ in $\A_\cC$ may be identified with the elements $[g] \in
G/CH$ such that $gHg^{-1}\leq H'$.

We have   canonical functors $\cC_G \to \A_\cC$ and $(\cC^{\op})_G \to
\bO_\cC^{\op}$, which in the above picture
take  $g$ to $[g]$. Note furthermore that we have a functor $\alpha: \A_\cC
\to \bO(G)^{\op}$ on objects given by $H
\mapsto G/CH$. In practice, our functors on $\A_\cC$ will be constructed in
this way from functors on $\bO(G)^{\op}$.
\enddemo

\numbereddemo{Definition}\label{bredoncohom}
Let $X$ be a $G$-space. We define a $G$-local coefficient system on $X$ to
be a functor $\cF: (\dio X)_G \to \Rmod$. Define a $G$-action on
$C^*(X;\cF)$,
 by letting $g$ act by $(gf)(\sigma) = \cF((1_{g\sigma},g))(f(g^{-1}))$.
Define the Bredon cohomology of $X$ with
coefficients in $\cF$ as $H_G^*(X;\cF) = H(C^*_G(X;\cF))$, where
$C^*_G(X;\cF)= C^*(X;\cF)^G$.
\enddemo

We may replace $C^*(X;\cF)$ by the corresponding normalized chain complex
for the purpose of calculations. If $M$ is an $RG$-module, then we denote
cohomology with values in the `constant' coefficient system $\sigma \mapsto
M$ by $H^*_G(X;M)$.

Note also that to any $G$-local coefficient system $\cF$ we may construct a
new coefficient system $\tilde \cF$ given by $\sigma \mapsto
\cF(\sigma)^{G_\sigma}$, where $G_\sigma$ denotes the isotropy subgroup of
$\sigma$. We have a natural transformation $\tilde \cF \to \cF$ which
satisfies $H^*_G(X;\tilde \cF) \pil{\simeq} H^*_G(X;\cF)$, and  \pagebreak
$\tilde{\tilde \cF} \pil{=} \tilde \cF$. (E.g., the constant coefficient
system $M$ corresponds to the generic coefficient system $H^0(-;M)$ in this
manner.)
Hence, for the purpose of Bredon cohomology, one could (as in
\cite{bredon67}) without loss of generality have restricted attention to
$G$-local coefficient systems which satisfied $\cF = \tilde \cF$.
However, in practice making this restriction seems slightly unnatural and
moreover the more general definition corresponds to the generality needed
when studying homology representations
(see \cite{RS85}, \cite[Ch.\ 7]{benson98}).
\vfill

 A {\it generic coefficient system} is a functor $\fF: {\bO}(G)^{\op} \to
\Rmod$. Note that for any $G$-space $X$ we have a canonical (covariant)
functor $\theta: (\dio X)_G \to  {\bO}(G)^{\op}$.
On objects this is given by $\sigma \mapsto G_\sigma$, and on morphisms it
is given by $(f :g\sigma \to \sigma',g) \mapsto (gG_{\sigma'}g^{-1} \geq
G_{\sigma},g)$, in the
setting of Remark \ref{grothrem}. A generic coefficient system induces
functorially a $G$-local coefficient system
$\cF$ on any $G$-space $X$ when $\cF = \fF \theta$. (In much of the
literature generic coefficient systems are
just called coefficient systems, or Bredon coefficient systems.) We call a
$G$-local coefficient system, which is
induced by some generic coefficient system, an {\it isotropy $G$-local
coefficient system}. (Bredon calls this a
{\it simple coefficient system} on $X$.) (Note that any functor $\fF:
\bO_\cC^{\op} \to \Rmod$ can be extended
(nonuniquely) to a functor $\fF': \bO(G)^{\op} \to \Rmod$ (see
\cite[X3.4]{maclane71}), so the exact domain of
definition is not a major issue.) When there is no danger of confusion we
will not distinguish between a generic
coefficient system and the isotropy $G$-local coefficient system it induces.
However, to avoid confusion we have
chosen not to give meaning to the term coefficient system without any
prefix.
 
In general, cohomology with values in a $G$-local coefficient system on $X$
is of course not a $G$-homotopy invariant of $X$. However, $G$-equivariant
Bredon cohomology with values in a generic coefficient system does define a
$G$-homotopy invariant cohomology theory on $G$-spaces (see e.g.\ 
\cite{bredon67}).
 \vfill
Recall that a Mackey functor is a generic coefficient system, but defined on
a category with more morphisms than the orbit category, where one also has a
`transfer' map (cf.\ e.g., \cite{may96}, \cite{webb99}). We actually only
directly need the Mackey functor $H^n(-;M)$ where the transfer map is the
usual transfer map in group cohomology. This is a so-called cohomological
Mackey functor, meaning that restriction from one subgroup to a smaller
subgroup followed by the transfer back equals multiplication by the index.

Given a $G$-space $X$ with a $G$-local coefficient system $\cF$, we define
Borel equivariant cohomology as $H^*_{hG}(X;\cF) = H^*_G(EG \times X;\cF)$,
where $EG \times X$ is given the diagonal $G$-action and coefficient system
pulled back via the projection map $EG \times X \to X$.

Given a $G$-category $\D$ and a functor $F: \D_G \to \Rmod$ (or $F:
(\D^{\op})_G \to \Rmod$), we naturally get a $G$-local coefficient system on
$|\D|$ via $\dio \D \to \D$ (or $\dio \D \to \D^{\op}$). \pagebreak

\numbereddemo{{R}emark} \label{padicremark} In most of the paper we work
over $\Z_{(p)}$, for the reason that we need multiplication by the index of
a Sylow $p$-subgroup to be invertible in $R$. (Any result stated over
$\Z_{(p)}$ holds verbatim over any $\Z_{(p)}$-algebra.) As a technical point
however, when dealing with the Steinberg complex, we need to work over
$\Z_p$ instead, in order to appeal to Webb's Theorem~\ref{webbthm}.
One can translate results back to $\Z_{(p)}$ in the following way: $\Z_p$ is
faithfully flat over $\Z_{(p)}$, so that  a
map of $\Z_{(p)}$-modules is an isomorphism if and only if it is an
isomorphism after tensoring with $\Z_p$.
Furthermore if $X$ is of finite type, then one checks that $H^*_G(X;\cF)
\otimes_{\Z_{(p)}}\Z_p \cong H^*_G(X;\cF
\otimes_{\Z_{(p)}}\Z_p)$, and the same statement is of course true
non-equivariantly. Hence, when showing certain
maps are isomorphisms, we can without problems work over $\Z_p$.
\enddemo
 
\demo{{\rm 2.5.} Higher limits}
Let $\Rm{\D}$ denote the category of functors $\D \to \Rmod$. By $\lim^*_\D$
we mean the right derived functors of the inverse limit functor $\lim_\D :
\Rm{\D} \to \Rmod$. We have the following fundamental identification.
\enddemo

\advance\theoremcount by 1
\proclaimtitle{cf.\ \cite[3.3]{GZ67}}
\proclaim{Proposition}\label{prop1} Let $F: \D \to \Rmod$ be any functor{\rm
.} Then
$$\lim_\D{}^* F \cong \Ext^*_{R\D}(R,F) \cong  H^*(\D;F).$$
\endproclaim 

\demo{Proof}  
First note that we can write $\lim^0_\D F = \Hom_{R\D}(R,F)$ so $\lim^*_\D F
= \Ext^*_{R\D}(R,F)$. Let $\D\downarrow d$ denote the over-category of $d
\in \D$ and consider the functor $C_*(|\D\downarrow-|;R): \D \to R \mbox{-
(chain complexes)}$, which to each $d \in \D$ associates the simplicial
chain complex of $|\D \downarrow d|$.
One easily checks that the degree $n$ of the cochain complex
$\Hom_{R\D}(C_*(|\D\downarrow-|;R),F)$ equals
$\prod_{x_0 \to \cdots \to x_n} F(x_n)$ with boundary map coming from the
simplicial structure on $|\D|$. In
particular we see that $\Hom_{R\D}(C_*(|\D\downarrow-|;R),-)$ preserves
exactness, so that
$C_n(|\D\downarrow-|;R)$ is projective in $\Rm{\D}$. But since $\D\downarrow
d$ has an initial object, and therefore
has contractible nerve,
$C_*(|\D\downarrow-|;R)$ defines a projective resolution of the constant
functor $R$ in $\Rm{\D}$.  This shows
$\lim^*_\D F$ equals $H^*(\D;F)$.
\enddemo 

\numbereddemo{{R}emark} \label{sheafremark} The cohomology $H^*(\D;F)$ can
be viewed as sheaf cohomology of a certain sheaf $\Phi$ on $\Top(|\D|)$, the
topological space associated to $|\D|$, constructed from $F$. In this remark
we sketch this construction. Suppose for concreteness that $F$ is
contravariant.
For each open set $U \subseteq  \Top(|\D|)$ let
$$\Phi(U) = \{  r  \in \prod_{\sigma \in |\D|, U \cap \sigma \neq
 \emptyset}F(\sigma) | F(\phi)(r_{\sigma'})=r_{\sigma} \mbox{ for all
$\phi:\sigma \to \sigma'$ in $\di\D$}\} .
$$
It is easy to see that this defines a sheaf on $\Top(|\D|)$, and if $y \in
\Top(|\D|)$, then the stalk at $y$ is equal to $F(\sigma)$, where $\sigma$
is the unique nondegenerate simplex such
that
$y$ lies in the open cell associated to $\sigma$. (An open zero cell is a
point.)  
 \enddemo

We now sketch how to see that $H^*(\D;F) = H^*(\Top(|\D|);\Phi)$.
We will restrict ourselves to the case where $\D$ is a poset---the general
case can be obtained (rather messily) by
passing to a double subdivision. If $\D$ is a poset $\Top(|\D|)$ is an
ordered simplicial complex. For $\sigma \in |\D|$,
let $\nstar(\sigma) \subseteq \Top(|\D|)$ be the open star of $\sigma$,
i.e., the union of all open cells whose boundary
contains~$\sigma$. Note that $\nstar(\sigma)$ is open and contractible, with
$\Phi(\nstar(\sigma)) = F(\sigma)$.
Furthermore the \v Cech complex associated to the open covering given by
$\nstar(\sigma)$, $\sigma \in |\D|$, equals
the cochain complex considered in the proof of Proposition~\ref{prop1}.
Since the \v Cech complex calculates sheaf
cohomology \cite{godement58}, we are done.

\demo{{\rm 2.8.} Some auxiliary $G$-categories}
Let $\iota: \bO_\cC \to \bO_{\cC \cup \{e\}}$ be the inclusion functor, and
consider the undercategory $E\bO_\cC = G/e\downarrow \iota$ (see \cite[p.\
46]{maclane71}).
Note that $E\bO_\cC$ has a $G$-action in the obvious way, on objects 
$g \cdot f = gf$. (Alternatively, we may view
$E\bO_\cC$ as the category with objects pairs $(H,x)$ where $H \in \cC$ and
$x \in G/H$, and with a unique morphism
$(H,x) \to (H',x')$ if there exists a $G$-map $G/H \to G/H'$ sending $x$ to
$x'$. The $G$-action in this picture is given
by $g\cdot (H,x) = (H,gx)$.) Likewise, set $E\A_\cC = \iota \downarrow G$,
where $\iota: \A_\cC \to \A_{\cC \cup
\{G\}}$ is the inclusion. We define a $G$ action via $g \cdot i = c_g i$,
where $(i: H \to G) \in E\A_\cC$ and $c_g$
denotes conjugation by $g$. (Note that the category $E\bO_\cC$ is isomorphic
to the category denoted ${\bf
X}^\beta_\cC$ in
\cite{dwyer98}, and $E\A_\cC$ is equivariantly equivalent to the category
denoted  ${\bf X}^\alpha_\cC$ in
\cite{dwyer98}.)
\enddemo

\advance\theoremcount by 1
\numbereddemo{{R}emark}
For any collection $\cC$ we have $G$-equivariant functors
$$ E\bO_\cC \to \cC \leftarrow E\A_\cC$$
given by $(f: G \to G/H) \mapsto G_f$ and $(i: H \to G) \mapsto i(H)$
respectively. These are (in general non-equivariant) equivalences of
categories with maps the other way being given by $H \mapsto (\mbox{proj}:G
\to G/H)$ and $H \mapsto (\mbox{incl}: H \to G)$ respectively. (This was a
key observation made by Dwyer
\cite{dwyer97}, \cite{dwyer98}.)
\enddemo

Note that any functor $F: \bO_\cC \to \Rmod$ induces a $G$-local coefficient
system on $|E\bO_\cC|$ via $(\dio E\bO_\cC)_G \to (\dio\cC)_G \to
(\cC^{\op})_G \to \bO_\cC^{\op}$, and is hence a pullback of the coefficient
system on $|\cC|$. This coincides with the isotropy $G$-local coefficient
system induced by $F$ viewed as a generic coefficient system, and is on
objects given by $(G\pil{f}G/H_0 \to \cdots \to G/H_n) \mapsto F(G_f)$. We
keep the notation $F$ for this $G$-local coefficient system on $|E\bO_\cC|$.
Likewise we get a $G$-local coefficient system on $|E\A_\cC|$ via $(\dio
E\A_\cC)_G \to (\dio\cC)_G \to \cC_G \to \A_\cC$, on objects given by $(H_0
\to \cdots H_n \pil{i} G) \mapsto F(i(H_n))$, and we also keep the notation
$F$ for this coefficient system. (Note that if $F$ is of
 the form $F = \fF \alpha$, where $\alpha$ is the canonical functor $\alpha:
\A_\cC \to \bO_\cC^{\op}$ and $\fF:
\bO(G)^{\op} \to \Rmod$ is some functor, then this is the isotropy $G$-local
coefficient system  induced by~$\fF$.)

\proclaim{Proposition}\label{prop2}
There are  natural isomorphisms $|E\bO_\cC|/G = |\bO_\cC|$ and  $|E\A_\cC|/G
= |\A_\cC|${\rm .}
 Furthermore $H^*(\bO_\cC;F)
= H^*_G(|E\bO_\cC|;F)$ and\break  $H^*(\A_\cC;F) = H^*_G(|E\A_\cC|;F).$
\endproclaim 

\vglue6pt

\demo{Proof} The first two statements follow  by inspection of the
simplices. But by the definition of the coefficient
systems and the first two statements, we now get that $C^*(|E\bO_\cC|;F)^G =
C^*(\bO_\cC;F)$ and $C^*(|E\A_\cC|;F)^G
= C^*(\A_\cC;F)$. Hence the last two statements follow by the definition of
Bredon cohomology.
\enddemo 

\section{Higher limits over $\cC$-orbit categories} \label{orbitsec}

In this section we examine higher limits over $\cC$-orbit categories.

\proclaim{Theorem}\label{thm1'}
Let $\cC$ be an arbitrary collection{\rm .} Let $F: \bO_\cC^{\op} \to \Rmod$
be a functor concentrated on conjugates of a fixed subgroup $H \in \cC${\rm
.} 
Then
\begin{itemize}
\ritem{1.} $H^*_G(|E\bO_\cC|;F) \cong H^*_{hWH}(|\cC_{\geq
H}|,|\cC_{>H}|;F(H)) \cong \tilde
H^*_{hWH}(\Sigma |\cC_{>H}|;F(H)) .$
\ritem{2.}
If $R = \Z_{(p)}$ and $|\cC_{>H}|^Q$ is $\F_p$-acyclic for all nontrivial
$p$-subgroups $Q$ in $WH$, then
\begin{eqnarray*}
H^*_G(|E\bO_\cC|;F) &\cong& H^*_G(|\cC|;\cF) \cong H^*_{WH}(|\cC_{\geq
H}|,|\cC_{>H}|;F(H))\\
& \cong &\tilde
H^*_{WH}(\Sigma |\cC_{>H}|;F(H)) ,
\end{eqnarray*}
\end{itemize}
where $\cF$ is the coefficient system obtained from $F$ by precomposing
with\break $(\dio \cC)_G \to (\cC^{\op})_G
\to
\bO_\cC^{\op}$.
\endproclaim 
 
\demo{Proof} 
We first perform some preliminary rewritings. By definition
$$C^n_G(|\cC|;\cF) = \left( \prod_{H_0 < H_1 < \cdots < H_n}F(H_0)\right)^G
.$$ Since $F$ is atomic this equals $(\prod_{H < H_1 < \cdots <
H_n}F(H))^{NH}$ by Frobenius reciprocity,
so that 
$$C^n_G(|\cC|;\cF) \cong C^n_{NH}(|\cC_{\geq H}|,|\cC_{>H}|;F(H)).$$
We want to do a similar rewriting with $E\bO_\cC$ in place of $|\cC|$.
View for a moment $E\bO_\cC$ as being replaced by its ($G$-equivalent)
skeletal subcategory with objects $G$-maps $G \to G/H'$ where $H'$ runs
through a set of conjugacy classes of objects in $\cC$. To simplify the
notation, set $Y = \cup_{H' \in \cC_{>H}}|E\bO_\cC^{H'}|$. As in the first
rewriting we have
\begin{eqnarray*}
 C^n_G(|E\bO_\cC|;F) &=& \left(\prod_{G \pil{f} G/H \to G/H_1 \to \cdots \to
G/H_n}F(G_f)\right)^G \\
                 &\cong& \left(\prod_{G \pil{f} G/H \to G/H_1 \to \cdots \to
G/H_n,  f(e) \in WH }F(H)\right)^{NH}\\
                 &=& C^n_{NH}(|E\bO_\cC^H|,Y;F(H)).
\end{eqnarray*}
These rewritings are natural, and give us corresponding identifications:
$$\begin{array}{ccc}
H^*_G(|\cC|;\cF)&\hskip.5in {\scriptstyle\cong}\hskip.5in  &
H^*_{WH}(|\cC_{\geq H}|,|\cC_{>H}|;F(H)) \\[4pt]
\Big\downarrow&&\Big\downarrow\\[9pt]
H^*_G(|E\bO_\cC|;F) &{\scriptstyle\cong} &  H^*_{WH}(|E\bO_\cC^H|,Y;F(H))\;
.\end{array}
$$ 
Note that on the right-hand side the coefficient systems are constant
coefficient systems, and hence come from the generic coefficient system
$H^0(-;F(H))$. 

We now prove $1$ of Theorem 3.1. First note that the $WH$ acts freely on the
pair $(|E\bO_\cC^H|,Y)$, or said
differently, the fixed-point pair under all nontrivial subgroups of $WH$
becomes a trivial pair of the form $(X,X)$.
Hence 
$$H^*_{WH}(|E\bO_\cC^H|,Y;F(H)) \cong H^*_{hWH}(|E\bO_\cC^H|,Y;F(H)) $$ by the
defining property of Borel cohomology. This
shows the first isomorphism in $1$, since the $WH$-map $(|E\bO_\cC^H|,Y) \to
(|\cC_{\geq H}|,|\cC_{>H}|)$ is a homotopy
equivalence of pairs. The second isomorphism follows directly from our
definition of the suspension.

We now assume $R = \Z_{(p)}$ and proceed to prove $2$. Consider again the
homotopy equivalence
\begin{equation}\label{he}
(|E\bO_\cC^H|,Y) \to (|\cC_{\geq H}|,|\cC_{>H}|).
\end{equation} 
Let $P$ be a Sylow $p$-subgroup of $WH$. Since $H^0(-;F(H))$ is a
cohomological Mackey functor we will be done, by an application of the
transfer, if we can see that this is a $P$-$\Z_{(p)}$-equivalence if and
only if $|\cC_{>H}|^Q$ is $\F_p$-acyclic for all nontrivial $p$-subgroups
$Q$ in $WH$. To this end, we have already observed that the $Q$-fixed point
set of the left-hand side of (\ref{he}) becomes the trivial pair. Also
$|\cC_{\geq H}|^Q$ is contractible since $\cC_{\geq H}^Q$ has the minimal
element $H$. Hence (\ref{he}) is a $P$-$\Z_{(p)}$-equivalence if and only if
$|\cC_{> H}|^Q$
is $\Z_{(p)}$-acyclic for all nontrivial\break $p$-subgroups $Q$ in $WH$.
Since $|\cC_{>H}|^Q$ is of finite type, we may
replace $\Z_{(p)}$ by~$\F_p$.
\enddemo 
\numbereddemo{{R}emark} By Smith theory (cf.\  e.g.\  \cite{may96}) it is in
Theorem~\ref{thm1'} enough to verify that $|\cC_{>H}|^Q$ is acyclic for all
subgroups $Q$ of order $p$.
\enddemo
\numbereddemo{{R}emark}\label{usefulrem}
Note that Theorem~\ref{thm1'} shows that if $F$ is concentrated on
conjugates of a subgroup $H$, then the equality $H^*_G(|E\bO_\cC|;F) \cong
H^*_G(|\cC|;\cF)$ is equivalent to the equality $\tilde H^*_{WH}(\Sigma
|\cC_{>H}|;F(H)) \cong \tilde H^*_{hWH}(\Sigma |\cC_{>H}|;F(H))$. The
condition in part $2$ of Theorem~\ref{thm1'} just gives an obvious criterion
for when relative Bredon and Borel cohomology with values in a cohomological
Mackey functor over $\Z_{(p)}$ coincide, namely when the action of a Sylow
$p$-subgroup is free (on the appropriate {\it pair}).
Alternatively the condition can be viewed as ensuring the collapse of the
$E_2$-term, a certain (relative) isotropy
spectral sequence, onto the {\it horizontal} axis (see \cite{dwyer98}).
Conditions of this type feature prominently in
the work of Webb \cite{webb87arcata}, \cite{webb87comment}, \cite{webb91},
and were also exploited by Dwyer
\cite{dwyer98}.
\enddemo

\numbereddemo{{R}emark} If $H$ is a maximal element in $\cC$, and $F$ is
concentrated on conjugates of $H$, then the formula of Theorem~\ref{thm1'}
shows that 
$$H^*_G(|E\bO_\cC|;F) \cong H^*(WH;F(H)).$$
In particular if $\cC$ contains a maximal element $H$ which is a non-Sylow
$p$-group, then $p\vert |WH|$, by a standard property of $p$-groups. Thus,
for example, when $F(H) = \Z_{(p)}$ with trivial $WH$-action,
$\lim^*_{\bO_\cC}F \cong  H^*(WH;\Z_{(p)})$ is nonzero in infinitely many
degrees. This shows that one does not in general get good vanishing results
for $\lim^*_{\bO_\cC}F$, when $\cC = \pA$, the collection of nontrivial
elementary abelian $p$-groups.
\enddemo

\proclaim{{C}orollary} \label{cor1'}
Let $F: \bO_\cC^{\op} \to \Zpmod$ be a functor, and assume that
$|\cC_{>H}|^Q$ is $\F_p$-acyclic for all $H \in \cC$
and for all nontrivial $p$\/{\rm -}\/subgroups $Q$ of $WH${\rm .}
 Then $H^*_G(|\cC|;\cF) \pil{\simeq} H^*_G(|E\bO_\cC|;F)${\rm .}
\endproclaim 

\demo{Proof}  Theorem~\ref{thm1'} establishes this corollary for atomic
functors. We want to establish the corollary for an arbitrary
functor $F$ by induction. Choose a maximal subgroup $H$ in $\cC$ on which
$F$ is nonzero, and let $F'$ denote the
subfunctor of $F$ obtained by setting $F'(H')=0$ if $H'$ is conjugate to $H$
and $F'(H') = F(H')$ otherwise. We have the
following diagram with exact rows:
\figin{2}{1000}
(The sequences remain exact after we take $G$-fixed points, since they are
split as modules.) 
 By induction on the number of conjugacy classes of subgroups in $\cC/G$ on
which the functor is nonzero, we can assume
\pagebreak that the statement is true for $F'$, and the statement is also
true for $F/F'$, since it is an atomic functor.
Hence the right and the left vertical arrows induce isomorphisms on
homology; so by the five-lemma this is also true
for the middle arrow.
\enddemo  

Before the next 
proof, recall that for a poset ${\cal P}$, and $x \in {\cal P}$, $\sta_{\cal
P} x$\break denotes the subposet $\{y \in
{\cal P}| x \leq y \mbox{ or } x \geq y\}$ and likewise $\link_{\cal P} x
=\break \{ y \in {\cal P} | x > y \mbox{ or } y <
x\}$.

\demo{Proof of Theorem~{\rm \ref{thm1}}}  We first prove the first part of
the theorem under the assumption that
$\cC$ is actually closed under passage to all $p$-overgroups, and will later
show that non-$p$-radical subgroups do
not play a role. We already know that $\lim^*_{\bO_\cC} F =
H^*_G(|E\bO_\cC|;F)$ from Propositions \ref{prop1} and
\ref{prop2}. Hence to prove the first part of the theorem, we just have to
show that the conditions of Corollary
\ref{cor1'} apply. By the assumption that $\cC$ is closed under passage to
$p$-overgroups, we have that for all $P \in
\cC$, $\cC_{>P} ={ \cal S}_p(G)_{>P}$. But for all nontrivial $p$-subgroups
$Q$ of $WP$,  ${\cal S}_p(G)_{>P}^Q $ has
contractible nerve. To see this note that for all $H \in {\cal
S}_p(G)_{>P}^Q$, identifying $Q$ with its preimage in
$NP$, $QH \in {\cal S}_p(G)_{>P}^Q$ so the inequalities $H \leq QH \geq Q$
provide a contracting homotopy (cf.\ also
\cite[1.5]{quillen78},\cite[Lemma 6.4.5]{benson98}). Hence Corollary
\ref{cor1'} shows that $H^*_G(|E\bO_\cC|;F)
\cong H^*_G(|\cC|;\cF)$.

Skipping slightly in the order in which the theorem was stated, we now
proceed to prove the atomic case. From
Theorem~\ref{thm1'} we know that $\lim^*_{\bO_\cC} F \cong \tilde
H^*_{NP}(\Sigma |\cC_{>P}|;F(P))$. By assumption
$\cC_{>P} = \pS_{>P}$. But the posets $\pS_{>P}$ and ${\cal S}_p(WP)$ have
$NP$-homotopy equivalent nerves, induced
by the $NP$-equiva\-riant maps $H \mapsto (H\cap NP)/P$ and $H/P \mapsto HP$ 
(cf.\ \cite[Prop.\
6.1]{quillen78}, \cite{benson98}), so
we can use the space $|{\cal S}_p(WP)|$ instead of $|\cC_{>P}|$. Writing
the definition of Bredon cohomology we now
obtain the wanted formula.

The last step in the proof is to see that we can replace the collection
$\cC$ (which we had assumed to be closed under passage to all
$p$-overgroups) by the smaller collection $\cC'$ as in the statement of the
theorem, without changing higher limits or Bredon cohomology. By filtering
the functor as in the proof of Corollary \ref{cor1'},
we see that it is enough to prove this when $F$ is an atomic functor,
concentrated on conjugates of a subgroup $P \in
\cC$.  If
$P
\in \cC \setminus \cC'$ then ${\cal S}_p(WP)$ is $WP$-contractible, a
contracting homotopy being given by $H \leq
HO_p(WP) \geq O_p(WP)$. Hence $H^*_G(|\cC|;F) =0$, by the formula for atomic
functors, which proves the wanted
result, since $F$ in this case is zero on $\cC'$. So, suppose that $P \in
\cC'$.  By Theorem~\ref{thm1'} it is enough to
prove that $\cC_{>P}$ and  $\cC'_{>P}$ are $WP$-homotopy equivalent. We will
do this by showing that we can
successively, in order of increasing size, remove $WP$-conjugacy classes of
subgroups in $\cC_{>P} \setminus
\cC'_{>P}$ from 
$\cC_{>P}$, without changing its $WP$-homotopy type.
Let $\tilde \cC$ be a $WP$-subposet of $\cC_{>P}$, with $\cC'_{>P} \subseteq
\tilde \cC \subseteq \cC_{>P}$. Suppose that $P' \in \tilde \cC \setminus
\cC'_{>P}$ and that $\tilde \cC_{>P'} = \cC_{>P'}$. Let $\hat \cC$ be the
poset obtained from $\tilde \cC$ by removing $WP$-conjugates of $P'$. It is
easy to check that we have a pushout square, which is a homotopy pushout of
 $NP$-space.
$$
\begin{array}{ccc}
 NP \times_{NP\cap NP'} |\link_{\tilde \cC}P'| &\lrar & |\hat \cC|  \\[4pt]
\Big\downarrow&&\Big\downarrow\\[9pt]
 NP \times_{NP\cap NP'} |\sta_{\tilde \cC}P'| &\lrar & |\tilde  \cC|.
\end{array}$$ 
Obviously $|\sta_{\tilde \cC}P'|$ is $NP'$-contractible.
Note that $\link_{\tilde \cC} P' = \cC_{>P'} \star \tilde \cC_{<P'}$, where
$\star$ denotes the join of posets (corresponding to the join of spaces; see
\cite{quillen78}). By the previous manipulations of posets, $|\cC_{>P'}|$ is
$NP'$ contractible since $P' \in \cC \setminus \cC'$. Hence  $|\link_{\tilde
\cC} P'|$ is $NP'$-contractible as well, so the above pushout square shows
that $|\hat \cC| \to |{\tilde \cC}|$ is
an $NP$-homotopy equivalence.

An induction now shows that $|\cC_{>P}|$ and $|\cC'_{>P}|$ are $NP$-homotopy
equivalent as wanted. (Compare also,
e.g., \cite[Prop 1.7]{TW91}, \cite[Thm.\ 7.1]{dwyer97}.)
\enddemo

\numbereddemo{{R}emark}In \cite{JMO92} the higher limits of a functor $F:
\bO_{\pSe} \to \Zpmod$ concentrated on the trivial subgroup $e$ with $F(e) =
M$ are denoted $\Lambda^*(G;M)$. So in this notation Theorem \ref{thm1} in
particular implies that
$$\Lambda^n(G;M) \cong H^{n-1}(\Hom_G(\St_*(G),M)).$$
\enddemo  
\numbereddemo{{R}emark} The $p$-radical subgroups have long been considered
in finite group theory, and e.g.\  appear in Gorenstein's $1968$ book, but
have more recently obtained a central role in connection with Alperin's
weight conjecture
 \cite{alperin87}, \cite{AF90}. Homotopy properties of the poset of
$p$-radical subgroups were first studied by Bouc
\cite{bouc84}. Their pertinence to the homotopy theory of classifying spaces
was (independently, and in the more
general context of compact Lie groups) discovered by
Jackowski-McClure-Oliver \cite{JMO92} who dubbed them
$p$-stubborn subgroups, since they just would not go away.
\enddemo

\numbereddemo{{R}emark} \label{pRequivrem}
Taking $\cC = \pSe$ and $P = e$ in the above proof, we see that $\pS$ and
$\pR$ are $G$-homotopy equivalent. This is a result of Bouc \cite{bouc84}
and 
Th\'evenaz-Webb \cite{TW91}.
\enddemo
\numbereddemo{{R}emark}\label{quillenrem1}
Note that Theorem~\ref{thm1} tells us that we can remove a {\it single}
conjugacy class $[H]$ from a collection $\cC \cap \pRe$ as in
Theorem~\ref{thm1} and still get the same result for {\it all} functors $F$
only if $\St_*(WH)$ is contractible. But a strong version of Quillen's
conjecture asserts that $O_p(W) = e$ implies that $\St_*(W)$ is
noncontractible,
so that minimality of the collection $\cC \cap \pRe$ in the sense above is
in fact equivalent to this version of
Quillen's conjecture. Note that, by Webb's Theorem \ref{webbthm}, $\St_*(W)$
is contractible if and only if $\tilde
H^*(|{\cal S}_p(W)|;\Z_{(p)})=~0$.  Quillen's conjecture, in its original
form states that $O_p(W)=e$ implies that $|{\cal
S}_p(W)|$ is non-contractible \cite[Conj. 2.9]{quillen78}. Quillen's
conjecture has been proven for a very large number
of groups in
\cite{AS93}, where they showed that $O_p(W)=e$ implies that $\tilde
H^*(|{\cal S}_p(W)|;\Z) \neq 0$ under some
restrictions on $W$ and $p$. (See also \cite[p.\ 206-207]{JMO92}.)
\enddemo
\numbereddemo{{R}emark}
The {\it non-equivariant} cohomology of subgroup complexes with values in a
$G$-local coefficient system has been studied previously
by, among others,  Ronan-Smith \cite{RS85}, \cite{RS86}, \cite{ronan89}
(cf.\ also \cite{benson98}), as a method for
constructing modular representations. It would be interesting to try to 
exploit the interplay between
the equivariant and non-equivariant cohomology.
\enddemo

\numbereddemo{{R}emark} \label{extendrem} The proof of Theorem~\ref{thm1}
reveals that the conclusion of the theorem holds not only for collections of
$p$-subgroups but also for all collections $\cC$
which are closed under extensions by $p$-groups. More precisely we need that 
if $H \triangleleft H' \leq G$ with $H \in \cC$ and $H'/H$ a $p$-group then 
$H' \in \cC$.
\enddemo

\demo{Proof of Theorem~{\rm \ref{pcentriccor}}}
The proof follows the pattern of the last part of the proof of
Theorem~\ref{thm1}, but uses the additional assumptions about $F$ together
with information about the structure of the Steinberg complex (for which we
refer to Section~\ref{st}) to remove the additional subgroups. As explained
in
Remark~\ref{padicremark} it is enough to prove the statement with $\Z_p$ in
place of $\Z_{(p)}$. Hence, for the rest
of the proof we adopt the convention that {\it all} coefficients are taken
over~$\Z_p$.

By filtering the functor as in the proof of Corollary \ref{cor1'}, it is
enough to prove the statement when $F$ is an atomic functor, concentrated on
conjugates of a subgroup $P \in \cC$.

Suppose first that $P \in \cC \setminus \cC'$. Since in this case
$H^*_G(|\cC'|;\cF)=0$ we need to see that
$H^*_G(|\cC|;\cF)=0$, i.e., by Theorem~\ref{thm1} that
$\Hom_{WP}(\St_*(W),F(P))$, or equivalently
$\Hom_{WP}(\SSt_*(W),F(P))$ (see Section~\ref{st}), is acyclic.

If $P$ is non-$p$-centric, then the kernel of the action of $WP$ on $F(P)$
has order divisible by $p$  and so, by
Corollary~\ref{stcritcor}, $\Hom_{WP}(\SSt_m(WP),F(P)) = 0$ for all $m$ and
$H^*_G(|\cC|;\cF) =0$ as wanted.

If $P$ is $p$-centric then $ZP$ is a Sylow $p$-subgroup in $CP$ and we can
write  $CP = ZP \times K$ where $p \hbox{$\not|$}\, |K|$. By Proposition
\ref{coinvprop} we have
$\St_*(WP)_K = \St_*(NP/PC(P))$, where $(\cdot)_K$ denotes coinvariants, and
hence $$\Hom_{WP}(\St_*(WP),
 F(P)) =
\Hom_{WP}(\St_*(NP/PCP),F(P)).$$
 But $O_p(NP/PCP) \neq e$ since $P \in \cC \setminus \cC'$  so that ${\cal
S}_p(NP/PC(P))$ is\break
$NP$-contractible and therefore  $\Hom_{WP}(\St_*(W),F(P))$ is acyclic.

Now suppose that  $P \in \cC'$.
It is, by Theorem~\ref{thm1'}, enough to prove that $$\tilde
H^*_{WP}(\Sigma|\cC_{>P}|;F(P)) = \tilde
H^*_{WP}(\Sigma|\cC'_{>P}|;F(P)),$$
 and that the same statement holds for Borel cohomology instead of Bredon
cohomology. As in the proof of Theorem~\ref{thm1}, we do this by showing
that we can successively, in order of
increasing size, remove $WP$-conjugacy classes of subgroups in $\cC_{>P}
\setminus \cC'_{>P}$ from
$\cC_{>P}$, without changing either    the Bredon or Borel equivariant
cohomology with coefficients $F(P)$.

Let the notation and assumptions be as in the last part of the proof of
Theorem~\ref{thm1}. The pushout square used there, together with the
Mayer-Vietoris sequence in Bredon cohomology, shows that we need to see that
$\tilde H^*_{NP\cap NP'}(\Sigma |\link_{\tilde \cC}P'|;F(P)) = 0$, and that
this also holds for Borel cohomology.
Now $\link_{\tilde \cC} P' = \cC_{>P'} \star \tilde \cC_{<P'}$. In
particular for all nontrivial $p$-subgroups $Q$ of $NP'$ we have that
$|(\link_{\tilde \cC} P')^Q|$ is contractible, since $|\cC_{>P'}^Q|$ is
contractible. Hence reduced Bredon and Borel
cohomologies coincide, and so we can restrict to considering Bredon
cohomology (see Remark \ref{usefulrem}).

By the definition of the join
$$\tilde C_*(|\link_{\tilde \cC} P'|) \cong \tilde C_*(|\cC_{>P'}|) \otimes
\Sigma \tilde C_*(|\tilde \cC_{<P'}|)$$
where $\Sigma$ denotes the suspension (i.e., dimension shifting) of chain
complexes. This chain complex is again $NP'$ chain homotopy equivalent to
$\SSt_*(WP') \otimes \Sigma\tilde C_*(|\tilde \cC_{<P'}|)$.
Assume that $P'$ is non-$p$-centric. Then we have that
\begin{eqnarray*} 
&&\hskip-.75in \Hom_{NP\cap NP'}(\SSt_m(WP') \otimes \tilde C_n(|\tilde
\cC_{<P'}|),F(P))\\[3pt]
& &\cong \
\Hom_{NP\cap NP'}(\SSt_m(WP'), \tilde C_n(|\tilde \cC_{<P'}|) \otimes F(P))
\\[3pt]
& &\cong \
\Hom_{NP'}(\SSt_m(WP'), (\tilde C_n(|\tilde \cC_{<P'}|) \otimes
F(P))\uparrow_{NP\cap NP'}^{NP'})\\[3pt]
& &= \ 0 .
\end{eqnarray*} 
Here the first isomorphism uses the fact
that permutation modules are self-dual and the last equality follows from
Corollary~\ref{stcritcor}, since $CP'$
acts trivially on $(\tilde C_n(|\cC_{<P'}|) \otimes F(P))\uparrow_{NP\cap
NP'}^{NP'}$ and that $P'$ is non-$p$-centric.
But this is true for all $m$, $n$, so we conclude that $\tilde H^*_{NP\cap
NP'}(\Sigma |\link_{\tilde \cC}P'|;F(P)) = 0$
as wanted.  

Now assume that $P'$ is $p$-centric and write $CP' = ZP' \times K$,
where\break $p  \hbox{$\not|$}\,  |K|$. We have
that
$|\link_{\tilde
\cC} P'|/K = |\cC_{>P'}|/K \star |\tilde \cC_{<P'}|$. But  $|\cC_{>P'}|/K$
is\break $NP'$-equivalent to $|{\cal
S}_p(NP'/P'C(P'))|$, by Proposition \ref{coinvprop}, which is
$NP'$-con\-tractible by the assumption that
$O_p(NP'/P'C(P')) \neq e$. Therefore
$$\tilde H^*_{NP\cap NP'}(\Sigma |\link_{\tilde \cC}P'|;F(P)) =  \tilde
H^*_{NP\cap NP'}(\Sigma |\link_{\tilde \cC}P'|/K;F(P)) = 0 ,$$
 also in this case, which finishes the proof of the theorem.
\enddemo

We note that the above proof can be simplified somewhat if one  is willing
to settle for removing {\it all}
non-$p$-centric subgroups, using the fact that the collection of $p$-centric
subgroups is closed under passage to
$p$-overgroups. \pagebreak

\numbereddemo{{R}emark} \label{quillenrem2} As in Remark \ref{quillenrem1}
the formula in Theorem \ref{pcentriccor} for the higher limits over atomic
functors
shows that minimality of the collection\break $\cC \cap \pD$ for all
functors which satisfy that $CP$ acts trivially
on
$F(P)$ is equivalent to a strong version of Quillen's conjecture.
\enddemo

\numbereddemo{Example} The $G$-homotopy type of $|\pD|$ is in general
different from that of $|\pS|$. For example if $G=\Sigma_3 \times \Z/2$ then
$|{\cal S}_2(G)|$ is $G$-contractible, whereas $|{\cal D}_2(G)| \simeq
G/(\Z/2 \times \Z/2)$.
\enddemo

\numbereddemo{{R}emark} \label{localconditions}For a given $p$-group $P$ and
a subgroup $Q \leq P$ we are interested in when we can have $ Q \in \pD$ for
some finite group $G$ with Sylow $p$-subgroup $P$. The definition of $\pD$
shows that necessary conditions include $$C_P(Q) =ZQ \mbox{ and } N_P(Q)/Q
\cap O_p(\Out(Q)) = e,$$ which are
conditions expressed just in terms of $P$ and $Q$.  Much more elaborate
conditions have been given in
\cite{puignotes}. Determining the possible principal $p$-radical $Q$'s along
with the possibilities for the groups
$N_G(Q)/QC_G(Q) \leq \Out(Q)$ is in fact equivalent to Brauer's $8$th
question \cite{brauer63} of determining the
possible fusion patterns in a $p$-group $P$. (See Section
\ref{stableelements} and
\cite{puig76}, \cite{puignotes}, \cite{thevenaz95}.)
\enddemo

\numbereddemo{Example} A concrete example, where $\pD$ is a proper
subcollection of the collection of $p$-radical
and $p$-centric subgroups, is given by letting $G$ equal the $p$-nilpotent
group $(\Z/2 \times \Z/2) \rtimes
\Sigma_3$, where $\Sigma_3$ acts on $\Z/2 \times \Z/2$ via $\Sigma_3 \to
\Z/2 \pil{{\rm flip}}
\Aut(\Z/2 \times \Z/2)$. In this case both the Sylow $2$-subgroups and $\Z/2
\times \Z/2$ are $2$-radical and
$2$-centric, but only the Sylow $2$-subgroups belong to ${\cal D}_2(G)$.
\enddemo

\demo{Proof of Theorem {\rm \ref{sscor}}}
Let $F_i$ be the subfunctor of $F$ obtained by setting $F(P) = 0$ for
$\hi(P)>i$, and let $\cF_i$ be the corresponding
$G$-local coefficient system. This induces an increasing filtration
$\{C^*(\cC';\cF_i)^G\}_i$ on $C^*(\cC';\cF)^G$. The
spectral sequence is now just the spectral sequence associated to a filtered
chain complex \cite{weibel94}. Since we
will need it in a while, let us write up explicitly the differential $d^1$
in the $E^1$-term:
$$C^k(|\cC'|;\cF_{i+1}/\cF_i)^G  \pil{d^1}
          C^{k+1}(|\cC'|;\cF_{i}/\cF_{i-1})^G .$$
More precisely, given
\begin{eqnarray*}
&&\alpha \in C^k(|\cC'|;\cF_{i+1}/\cF_i)^G\\[5pt]
&&\hskip.5in =
 \bigoplus_{[P'] \in \cC'/G;\hi(P')=i+1} \Hom_{NP'}(C_k(|\cC'_{\geq
P'}|,|\cC'_{> P'}|),F(P')),
\end{eqnarray*}
we get by conjugation a unique element (also denoted $\alpha$), which lives
in the same sum, but is now taken over all
$P'$, and not just conjugacy classes. We can describe $d^1(\alpha) $ as the
following composition:
\begin{eqnarray*}
d^1(\alpha): C_{k+1}(|\cC'_{\geq P}|,|\cC'_{>P}|) &\pil{\partial}&
\bigoplus_{P<P';\hi(P') = \hi(P)+1}
C_k(|\cC'_{\geq P'}|,|\cC'_{>P'}|)  \\
&\pil{\alpha}&\bigoplus_{P<P';\hi(P') = \hi(P)+1} F(P') \pil{F(P<P')} F(P)
\end{eqnarray*}
where $\partial$ is the boundary map associated to the triple $$(|\cC'_{\geq
P}|,
|\cC'_{\geq P,\hi\geq\hi(P)+1}|,|\cC'_{\geq P,\hi\geq\hi(P)+2}|).$$

That we can replace $\cC \cap \pRe$ by $\cC \cap \pD$ and $NP/P$ by
$NP/PC(P)$ if $CP$ acts trivially on $F(P)$ for all $P \in \cC$, follows
from Theorem~\ref{pcentriccor}.
\enddemo 

\section{The case of a finite group of Lie type} \label{liesec}

The purpose of this section is to examine the spectral sequence in
Theorem~\ref{sscor} in the special case of a finite
group of Lie type $G$ of characteristic $p$, where it degenerates into a
chain complex.

We first need to recall some standard Lie notation.
To a finite group of Lie type $G$ of characteristic $p$, we can associate a
split $BN$-pair of characteristic $p$ (cf.\ e.g.\  \cite[p.\ 50]{carter85},
\cite{CR2}), which we will always assume to be
fixed. Let $S$ be the set of simple reflections of the Weyl group of the
$BN$-pair. (Beware that this Weyl group need
not coincide with the Weyl group of the ambient algebraic group if the
Frobenius map is twisted.) Any parabolic
subgroup in $G$ is conjugate to exactly one `standard' parabolic subgroup,
and the standard parabolic subgroups in $G$
are in one-to-one correspondence with subsets $J \subseteq S$, and are
denoted $P_J$ correspondingly. One defines
$U_J$ as the largest normal $p$-subgroup in $P_J$. Recall that  $N_G(U_J) =
P_J$ and $P_J = U_J \rtimes
L_J$, where $L_J$ again is a finite group of Lie type, with a split
$BN$-pair of rank $|J|$. To any finite group of Lie
type $G$ of characteristic $p$ we can associate the Steinberg module $\St_G$
(see  \cite{humphreys87}, \cite{CR2}).
One of the many definitions of $\St_{G}$ is $\St_{G} =
H_{r-1}(\Delta;\Z_{(p)}) =  H_{r-1}(\Delta;\Z) \otimes
\Z_{(p)}$, where
$\Delta$ is the Tits building of the $BN$-pair of $G$ and $r$ is the rank of
the $BN$ pair. It is a projective
indecomposable module of dimension $|G|_p$, whose reduction mod $p$ is also
simple. 

A theorem of Borel and Tits \cite{BT71}, \cite{BW76} says that any
$p$-radical subgroup is conjugate to exactly one
$U_{J}$, and it easily follows (cf.\ \cite{TW91} and Remark \ref{btremark})
that $\pRe$ is equal to the opposite poset
of the poset of parabolic subgroups in $G$. Hence the $NP/P$ which occur in
Theorem \ref{sscor} will be of the form
$L_J$ and, e.g., by Webb's Theorem~\ref{webbthm},  we may replace
$\St_*(L_J)$ with $\St_{L_J}$ in degree $|J|-1$
(cf.\ Remark \ref{padicremark}). When the height function $\hi(U_I) = -|I|$
is chosen, Theorem \ref{sscor} now takes
the following simple form.

\proclaim{Theorem} \label{liecor}
Let $G$ be a finite group of Lie type of characteristic $p${\rm ,} $\cC$ be
a collection of $p$\/{\rm -}\/subgroups
 closed under passage to $p$\/{\rm -}\/radical overgroups{\rm ,} and $F:
\bO_\cC^{\op} \to \Zpmod$ any functor{\rm .}
 Then
$\lim_{\bO_\cC}^{*}F$ equals the homology of a cochain complex
$(C^*(F),\delta)$ with
$$C^i(F) =  \bigoplus_{J \subseteq S ; |J|=i ; U_J \in \cC }
\Hom_{L_J}(\St_{L_J},F(U_J)).$$
\endproclaim 

We now proceed to examine the boundary map in Theorem \ref{liecor}, which
will be the content of Theorem \ref{identprop}. It turns out that special
properties of finite groups of Lie type
make this boundary map very tractable.

Let $\Rm{{\bO}(G)}$ denote the category of functors $\bO(G)^{\op} \to
\Rmod$. 
For any subgroup $H$ of $G$, we have a functor $\Rm{{\bO}(G)} \to \Rm{WH}$
given by evaluation $F \mapsto F(H)$. This functor has a right adjoint given
by associating to a $WH$-module $M$ the functor $E_{H,M}$ defined as 
$$E_{H,M}(H') = (\oplus_{\mor_G(G/H,G/H')}M)^{WH}.$$
Here $WH$ acts on $M$ in the prescribed way, and on $\mor_G(G/H,G/H')$ via
$g \cdot \alpha = \alpha \circ g^{-1}$.
($\mor_G$ denotes morphisms in the category of left $G$-sets.) In
particular, for any functor $F$ we get a natural
transformation $F \to E_{H,F(H)}$ induced by the unit of this adjunction.
(See \cite{GM92} for a careful explanation of
the role of such adjunctions.)

For finite groups of Lie type, this right adjoint has a particularly simple
form.

\proclaim{Lemma} \label{oplemma}
 If $G$ is a finite group of Lie type of characteristic $p${\rm ,} and 
$U_I${\rm ,} $U_{J'}$ are
the unipotent radicals of standard
parabolic subgroups
$P_I$, $P_{I'}${\rm ,} then
$$ E_{U_I,M}(U_{I'}) = \tuborg 0 & \mbox{if } I' \not \subseteq I \\
                               M^{U_{I'}} &\mbox{if } I' \subseteq I.
\sluttuborg  
$$
\endproclaim

\demo{Proof} 
As usual we have
$$\mor_G(G/U_I,G/U_{I'})= \{ g \in G | g^{-1}U_Ig \leq U_{I'}\}/U_{I'}.$$
We claim that
\begin{equation}\label{btequation} \{ g \in G | g^{-1}U_Ig \leq
U_{I'}\}/U_{I'} = \tuborg P_I/U_{I'} & \mbox{if } I' \subseteq I \\
\emptyset
&\mbox{otherwise} \sluttuborg.
\end{equation}

To show this we follow standard $BN$-pair notation as in \cite[\S 64, 65,
69]{CR2} (except we use $U_I$ for their
$V_I$), and assume some familiarity with this material. If $g \in P_I$, it
is clear that $g^{-1}U_Ig \subseteq  U_{I'}$
if and only if  $I'
\subseteq I$, and so we  suppose $g \not \in P_I$. Using the Bruhat
decomposition \cite[Thm.\ 65.4]{CR2} we can write
$g = b\dot wb'$, where $\dot w$ is a lift of a $w \in W \setminus W_I$, and
$b,b' \in B$. Hence we can without
restriction assume $g =\dot w$, where $w$ has an expression as a reduced
word of simple reflections, starting with
$s_\alpha$, for some $\alpha \in S \setminus I$. But then $w^{-1}\alpha$ is
a negative root, by the theory of root
systems
\cite[Thm.\ 64.16(vii)]{CR2}. This means that $g^{-1}U_Ig \not \subseteq
U_{\emptyset}$ \cite[69.5]{CR2}, and the
claim is justified.

We now see that
\medbreak
\hfill ${\displaystyle E_{U_I,M}(U_{I'}) = (\oplus_{P_I/U_{I'}}M)^{P_I/U_I}
=
\Hom_{P_I}(\Z_{(p)},M\uparrow^{P_I}_{U_{I'}}) = M^{U_{I'}} .}$
\enddemo

\numbereddemo{{R}emark}\label{btremark}
The above lemma can be used to show that the poset of nontrivial $p$-radical
subgroups in $G$ is naturally
$G$-equivalent to the opposite poset of the poset of parabolic subgroups in
$G$. Namely, by the Borel-Tits theorem
\cite{BT71} there is a one-to-one correspondence between the objects being
given by taking normalizer and
$O_p(\cdot)$ respectively. Now (\ref{btequation}) and \cite[Thm.\
65.19]{CR2} combine to show that this actually is an
equivalence of {\it posets}. (This fact is of course well known, although a
complete proof does not seem to be in the
literature.)
\enddemo

\numbereddemo{{R}emark} The proof of Lemma \ref{oplemma} easily adapts to
give a proof of another related fact. For $G$ a finite group of Lie type and
$\cC$ the collection of $p$-subgroups which contain a conjugate subgroup of
$U_I$,
the category $\bO_{\cC \cap \pRe}$ is equivalent to $\bO_{{\cal
B}_p^e(L_I)}$. Note that $L_I$ has smaller
semisimple rank than $G$ if $I \neq S$, which allows for proofs by
induction. One can for instance use this to give a
different inductive proof of Theorem \ref{liecor} and its addendum
Theorem~\ref{identprop}, which instead of
subgroup complexes only relies on knowledge of the irreducible
representations of $G$.
\enddemo

\proclaim{Lemma} \label{stlemma}  $(\St_{\Z_{(p)}G})^{U_J} \cong
\St_{\Z_{(p)}L_J}$ as
 $\Z_{(p)}L_J$\/{\rm -}\/modules{\rm .}
\endproclaim

\demo{Proof} 
The corresponding statement over $\F_p$ instead of $\Z_{(p)}$ is a
consequence of a general theorem of S. Smith about $U_J$-fixed points of
irreducible modules of finite groups of Lie type (cf.\
\cite{smith82}, \cite{cabanes84}).
We now show how to get the statement over $\Z_{(p)}$ from that over $\F_p$,
in a slightly indirect way. Note that
$\Ext^1_{U_J}(\St_{\Z_{(p)}G},\Z_{(p)}) = 0$ since $\St_{\Z_{(p)}G}$ is
projective. Hence we have an exact sequence
$$ 0 \to \Hom_{U_J}(\St_{\Z_{(p)}G},\Z_{(p)}) \pil{p}
\Hom_{U_J}(\St_{\Z_{(p)}G},\Z_{(p)}) \to  \Hom_{U_J}(\St_{\Z_{(p)}G},\F_p)
\to 0$$
which shows that $(\St_{\Z_{(p)}G})^{U_J} \otimes_{\Z_{(p)}} \F_p \cong
(\St_{\F_pG})^{U_J}$, using the fact that
$\St_{\Z_{(p)}G}$ is self-dual.

Hence the two expressions in the lemma are isomorphic after reduction mod
$p$. Since $\St_{\Z_{(p)}L_J}$ is
$\Z_{(p)}L_J$-projective this isomorphism can be lifted over to
$\Z_{(p)}L_J$, where it has to be an
isomorphism by Nakayama's lemma.
\enddemo 
\pagebreak

\proclaim{Lemma} \label{bdlemma} Assume $J' \subset J${\rm ,} $|J| = |J'| +
1${\rm ,}
 and let $F = E_{U_J,\St_{L_J}}${\rm .} Then the restriction and projection
of the boundary map
  $\delta: C^{|J'|}(F) \to
C^{|J|}(F)$ onto the corresponding factors $J'$ and $J${\rm ,}
$$ \delta: \Hom_{L_{J'}}(\St_{L_{J'}}, (\St_{L_J})^{U_{J'}}) \to
\Hom_{L_{J}}(\St_{L_J}, \St_{L_J})$$
is an isomorphism{\rm .}
\endproclaim

\demo{Proof}  First note that by Lemma \ref{stlemma} we know that
$(\St_{L_J})^{U_{J'}} \cong \St_{L_{J'}}$.
By Nakayama's lemma it is enough to prove the stated isomorphism after
reducing both sides modulo $p$. Furthermore, one immediately
checks, as in Lemma \ref{stlemma}, that the appropriate $\Ext^1$ vanishes so
that we may commute
$-\otimes_{\Z_{(p)}}\F_p$ and $\Hom$. We have hence reduced the question of
showing the isomorphism  to a
question over $\F_p$, and we therefore make the notational convention that
in the rest of this proof, everything will
be taken over $\F_p$ instead of $\Z_{(p)}$ which we suppress   from the
notation.

Since $\St_{L_{J'}}$ is absolutely irreducible over $\F_p$,
 $\Hom_{L_{J'}}(\St_{L_{J'}}, \St_{L_{J'}}) = \F_p$ (cf.\ \cite[\S
13.4]{bourbakialgebreVIII}), and likewise for $J$. Hence to show the
isomorphism, it is enough to show that the identity map of $\St_{L_{J'}}$ is
carried to a nonzero map under $\delta$.

The proof of Theorem \ref{sscor} gives a concrete description of $\delta$.
Since the homology modules of the chain complexes appearing in the proof 
of Theorem
\ref{sscor} are projective, we can replace the chain complexes with their
homology. Hence 
the restriction and projection of the boundary map  $\delta: C^{|J'|}(F) \to
C^{|J|}(F)$ onto the corresponding factors
$J'$ and $J$ is given as follows:
$$\delta(1): \St_{L_J} \to  \St_{L_{J'}}\uparrow_{P_{J'}}^{P_J} \pil{1}
\St_{L_{J'}}\uparrow_{P_{J'}}^{P_J}
 \to \St_{L_J}.$$
The first map is injective, by the long exact sequence in homology inducing
it. The last map is surjective, since it is
seen to be a nontrivial map of $P_J$-modules by its definition.  But it is
well known that
$\St_{L_{J'}}\uparrow_{P_{J'}}^{P_J}$ contains a {\it unique} summand of
$\St_{L_J}$ (cf.\ \cite[Thm.\ 67.10]{CR2}),
and so we conclude that $\delta(1)$ is in fact nontrivial, and hence is an
isomorphism.
\enddemo 

\proclaimtitle{Addendum to Theorem \ref{liecor}}
\proclaim{Theorem} \label{identprop} \hskip-8pt
Assume $J' \subset J$, $|J| = |J'| +~1${\rm ,}
 and let $e$ be a central idempotent in $\Z_{(p)}L_J$ corresponding to the
Steinberg block {\rm \cite{CR2}.}
 The restriction
and projection of the boundary map $\delta: C^{|J'|}(F) \to  C^{|J|}(F)$
onto the corresponding factors $J'$ and $J$ are
given as the top horizontal composite in the following diagram\/{\rm :}
\figin{3}{875}%
\noindent Here the top left horizontal map is induced by the map in $\bO_{\cC'}$
{\rm (}\/landing in the invariants\/{\rm ),} and $\delta'$ and $\delta''$
are the boundary maps in the chain complex
for
$E_{U_J,F(U_J)}$ and
$E_{U_J,eF(U_J)}$ respectively{\rm .}
\endproclaim   

\demo{Proof} 
 The diagram is commutative, since it is induced by the natural
transformations of functors $F \to E_{U_J,F(U_J)} \to
E_{U_J,eF(U_J)}$. That $\delta''$ is an isomorphism  is the content  of
Lemma \ref{bdlemma}.
\enddemo 

\section{The Steinberg complex}\label{st}

\hglue-1pt In this section we analyze the Steinberg complex $\St_*(G)\! =\!
\tilde C_*(|\pS|;\Z_{(p)})$. Many of the
results in this section concerning the Steinberg complex are known to group
theorists, and  are mainly due to Webb
\cite{webb87arcata}, \cite{webb87comment}, \cite{webb91} (see also
\cite{benson98}, \cite{YW94}, \cite{bouc92a}, \cite{bouc92b},
\cite{bouc96}).  However we find it worthwhile to
include them, since they are important for calculations with the results of
the previous section.
Here we shall mainly work over the $p$-adic integers $\Z_p$. The main reason
for this is the following
theorem.

\proclaimtitle{Webb \cite{webb91}}
\proclaim{Theorem} \label{webbthm} The Steinberg complex over $\Z_p$ has the
following structure
$$\St_*(G;\Z_p) = P_* \oplus D_*$$
as $\Z_{p}G$\/{\rm -}\/modules{\rm ,} where $P_*$ is a complex of
$\Z_{p}G$\/{\rm -}\/projective modules{\rm ,}
and $D_*$ is a $\Z_{p}G$\/{\rm -}\/split acyclic complex{\rm . }
\endproclaim 

(Unfortunately, the proof given in \cite{webb91} is not very clear on how to
pass from the corresponding statement over $\F_p$, which is first proved, to
the statement over
$\Z_p$; the problem is that although all morphisms lift, the lifting is not
functorial---this was pointed out to me
by P. Symonds, who also explained   how to overcome it
\cite{symondsletter}.) In the category of
bounded chain complexes over $\Z_pG$ the Krull-Schmidt theorem holds, so
that there exists a unique minimal such
$P_*$ (unique up to {\it isomorphism} of chain complexes) (see
\cite[5.17.1]{benson98}). We denote such a minimal
$P_*$ by
$\SSt_*(G)$---it is arguably really $\SSt_*(G)$ which deserves the name the
Steinberg complex \cite{alperin90}.

We now proceed to give restrictions on which projective modules   can occur
in $\SSt_*(G)$. Since in any degree
$\St_*(G;\Z_p)$ is a sum of permutation modules, we first need a method of
seeing which projective modules can
occur as summands of a given permutation module.  Recall that there is a
bijection between projective modules over
$\Z_p$ and $\F_p$ (cf.\ \cite[14.4]{serre77}),  and so we can restrict our
attention to fields. Our method is the
following lemma (with $T = k$), which is hinted at in \cite[Prop.\
5.3]{webb87arcata}. (See also \cite{robinson89}.)

\proclaim{Lemma} \label{grplemma} Let $k$ be a field of characteristic
$p${\rm ,}
 and assume that $H$ is a subgroup of a finite group $G${\rm .} Let $S$ be a
simple $kG$\/{\rm -}\/module
 and let $T$ be a simple
$kH$-module{\rm .} Denote the projective covers $P_S$ and $P_T${\rm .} Then
$[\End_GS:k]$ times the multiplicity of $P_S$ in $T\uparrow^G_H$ equals
$[\End_HT:k]$ times the multiplicity of $P_{T}$ in $S\downarrow^G_H${\rm .}
\endproclaim

\demo{Proof}  Let $\Omega^0M$ denote the largest submodule of a $kG$-module
$M$, which does not contain any
projective summands (unique only up to isomorphism). There is a
decomposition  $M = \Omega^0M \oplus P$, where $P$
is projective. Note that the multiplicity of $P_T$ in $S\downarrow^G_H$ is
equal to the $\End_HT$-dimension of
$\coker(\Hom_H(T,\Omega^0(S\downarrow^G_H)) \hookrightarrow
\Hom_H(T,S\downarrow^G_H))$. Likewise the
multiplicity of $P_S$ in $T\uparrow_H^G$ equals the $\End_G(S)$-dimension of
$\coker(\Hom_G(\Omega^0(T\uparrow_H^G),S) \hookrightarrow
\Hom_G(T\uparrow_H^G,S))$. So in order to prove the
lemma, we just need to justify that $\Hom_H(T,\Omega^0(S\downarrow^G_H))
\cong
\Hom_G(\Omega^0(T\uparrow_H^G),S)$. We will show this by showing that we may
replace $\Hom$ in both expressions
with $\underline\Hom$, i.e., $\Hom$ in the stable module category, where the
statement will follow by Frobenius
reciprocity in the stable module category (cf.\ \cite[p.\ 74]{alperin86}).
Assume $\alpha: \Omega^0(k\uparrow_H^G)
\to S$ is a $kG$-map which factors through a projective module, and hence
through $P_S$.
\figin{4}{1000}%
\noindent If $\alpha$ is nontrivial then $\Omega^0(T\uparrow_H^G) \to P_S$ has to be
surjective, since it is
surjective on the head, contradicting that $\Omega^0(T\uparrow_H^G)$ does
not contain any projective summands.
Hence $\alpha$ is trivial and  $\Hom_G(\Omega^0(T\uparrow_H^G),S) \cong
{\underline\Hom}_G(T\uparrow_H^G,S)$.
The proof that
$\Hom_H(T,\Omega^0(S\downarrow^G_H)) \cong
{\underline\Hom}_H(T,S\downarrow^G_H)$ is dual.
\enddemo   

By \cite{webb91} the conclusion of Webb's Theorem~\ref{webbthm} does not
change if we replace $\St_*(G;\Z_p)$ by a $G$-homotopy equivalent chain
complex. Hence
$\SSt_*(G)$ has to be a direct summand of the reduced normalized cellular
chain complex of any space $X$ $G$-homotopy equivalent to $|\pS|$.  Lemma
\ref{grplemma} says that if the space $X$ has `large' isotropy groups, then
the simple modules $S$, corresponding to projective modules $P_S$ appearing
in $\SSt_*(G)$, also have to be `large'. Note that a nondegenerate
$n$-simplex $\sigma = (P_0 < \cdots < P_n) \in |\cC|_n$ has stabilizer
$G_\sigma = \cap_{i=0}^nNP_i$.

There are many spaces which are known to be $G$-homotopy equivalent to
$|\pS|$. This is true for both $|\pR|$
and $|\pA|$ (\cite{bouc84}, \cite{quillen78}, \cite{TW91}; for proofs see
also the proofs of Theorem~\ref{thm1} and
\ref{thm2}). If $G$ is a finite group of Lie type of characteristic $p$ then
$|\pR|$ is just the subdivision of the
building $\Delta$ of $G$ (see Remark \ref{btremark}), and in particular they
are $G$-homotopy equivalent.
Furthermore $|\pA|$ is $G$-homotopy equivalent to the smaller complex given
by the nerve of the collection of
elementary abelian $p$-subgroups $V$ satisfying that $V$ equals the elements of order $p$ in the
center of
$CV$ \cite[\S 6.6]{benson98}.  The so-called {\it commuting complex}
$K_p(G)$ of Alperin
\cite{alperin90}, \cite{aschbacher93}, which has $n$-simplices
$(n+1)$-tuples of pairwise commuting subgroups of order $p$ and is also
$G$-equivalent to $|\pS|$ \cite[p.\
187]{thevenaz93}. (One way to see this is to note that the category with
objects, arbitrary tuples of pairwise
commuting subgroups of order $p$ viewed as a poset under inclusion, is
equivariantly equivalent to $\pA$; functors
back and forth are given by associating to an $n$-tuple the elementary
abelian $p$-subgroup it generates, and to an
elementary abelian $p$-group associating the tuple of all rank one
$p$-subgroups.)

Fix a Sylow $p$-subgroup $P$ of $G$ and define the following subgroup:
$$ B = \cap_{Q\in \pR_{\leq P}} NQ.$$
By Sylow's theorem it is uniquely defined up to conjugacy. It lies between
$CP$ and $NP$ and so $p| |B|$ if $p | |G|$. If $G$ is a finite group of Lie 
type of characteristic $p$ then $B = NP$ and $B$ is a Borel subgroup.

\proclaim{Proposition} \label{stcrit}
Let $P_S$ be the projective module over $\Z_pG$ corresponding to the
projective cover over $\F_pG$ of a simple
$\F_pG$-module $S${\rm .} Assume $P_S$ appears as a summand in
$\SSt_m(G)${\rm .}
 Then we can say the following about $S${\rm :}
\begin{itemize}
\ritem{1.} No element of order $p$ in $G$ acts trivially on $S${\rm .}
\ritem{2.} The module $S\downarrow^G_B$ contains the projective cover
$P_{\F_p}$ of the trivial\break
$\F_pB$\/{\rm -}\/module
$\F_p$ as a direct summand. In particular $\dim_{\F_p}S \geq |B|_p${\rm .}
\ritem{3.} There exists an elementary abelian $p$\/{\rm -}\/subgroup $V$ of
$G$ with $\rk V \geq m+1$ such that
$S\downarrow^G_{CV}$ contains the projective cover $P_{\F_p}$ of the
trivial\break $CV$\/{\rm -}\/module $\F_p$
as a direct summand{\rm .} In particular  $\dim_{\F_p}S \geq |CV|_p \geq
p^{m+1}${\rm .}
\end{itemize}

\endproclaim  

\demo{Proof}  
To see the first claim let $K = \ker(G \to \Aut(S))$ and assume that $P_S$
is a summand of $\Z_{p}[G/G_\sigma]$ for some $\sigma \in |\pS|$.
Pick a Sylow $p$-subgroup $P$ containing the chain of $p$-subgroups defining
$\sigma$. Since $K$ is normal in $G$, $P \cap K$ is a Sylow $p$-subgroup of
$K$. 
Now suppose there exists an element of order $p$ which acts trivially on
$S$, so that $ P \cap K \neq e$. Then, since
$P$ and $P\cap K$ are $p$-groups, we can find a nontrivial element $g \in
P\cap K$ which is fixed under the
conjugation action of $P$. This means that $\langle g\rangle \subseteq
G_\sigma$. By Lemma \ref{grplemma},
$S\downarrow^G_{G_\sigma}$ contains a summand of the projective
$\F_pG_\sigma$-module $P_{\F_p}$. But then the
further restriction $S\downarrow^G_{\langle g\rangle}$ has to contain a
projective module as well. Hence it cannot
be the trivial module, which is a contradiction since $g \in K$.

To see the second claim we choose to view the Steinberg complex as the
reduced normalized chain complex of $|\pR|$.
If $P_S$ appears as a direct summand in $\tilde C_*(|\pR|;\Z_p)$, then $P_S$
appears as a direct summand of $\Z_{(p)}[G/G_\sigma]$, for some $\sigma \in
|\pR|$.
Hence by Lemma \ref{grplemma}, $S\downarrow^G_{G_\sigma}$ contains the
projective cover $P_{\F_p}$  of the trivial $\F_pG_\sigma$-module $\F_p$ as
a direct
summand; so by Sylow's theorem $G_\sigma$ contains a conjugate of $B$, and
$S\downarrow^G_B$ has to contain a
copy of the projective cover of the trivial $\F_pB$-module as well.  It is
now obvious that $\dim_{\F_p}S \geq |B|_p$,
since a projective $\F_pB$-module has dimension divisible by~$|B|_p$.

To see the last claim we choose to view the Steinberg complex as the reduced
normalized chain complex of $|\pA|$. Observe that if $\sigma = (V_0 < \cdots
< V_m)$ is a nondegenerate $m$-simplex in $|\pA|$ then $CV_m \leq G_\sigma$
and $V_m$ has at least rank $m+1$. Hence the last claim follows from Lemma
\ref{grplemma} as well.
\enddemo 

Note that $(1)$ in the above proposition, combined with Theorem~\ref{thm1},
has \cite[Prop.\ 5.5]{JMO92} as a consequence (compare also \cite[Thm.\
6.3]{dwyer97}). Let us for reference write a consequence of 
Proposition~\ref{stcrit}$(1)$ and $(3)$ out explicitly---the part coming from 
Proposition \ref{stcrit}$(3)$ is new.
\proclaim{{C}orollary} \label{stcritcor}
Let $G$ be a finite group and let $M$ be a $\Z_pG$\/{\rm -}\/module{\rm .}
Assume there exists a filtration 
$0 = N_0 \subseteq N_1 \subseteq \cdots \subseteq N_m = M$ of $M$ 
such that for each $i$ either 
\begin{itemize}
\ritem{1.} $\ker(G \to \Aut(N_i/N_{i-1}))$ has order divisible by $p${\rm ,}  
or
\ritem{2.} $N_i/N_{i-1}$ is generated over $\Z_p$ by strictly less that 
$p^k$ elements,
\end{itemize}
then $\Hom_G(\SSt_{k-1}(G),M)=0${\rm .}

In particular{\rm ,} by Theorem {\rm \ref{thm1},}
 if $F$ is a functor $F: \bO_{\pSe}^\op \to \Zpmod$ concentrated on the trivial
subgroup $e$ with $F(e)=M${\rm ,} where $M$ satisfies the above 
condition{\rm ,} then $\lim^k_{\bO_{\pSe}}F =0${\rm .}
\endproclaim

\numbereddemo{Example}
If $G$ is finite group of Lie type in characteristic $p$, then the Steinberg 
module $\St_G$ is the
{\it only} simple module $S$ such that $S\downarrow^G_B$ contains a summand
of $P_{\F_p}$, agreeing nicely with our knowledge that the Steinberg complex
in this case is chain homotopy equivalent to the Steinberg module
concentrated in degree one less than the semi-simple rank of $G$.
\enddemo

\numbereddemo{Example} The structure of the Steinberg complex for the
symmetric group $\Sigma_n$ is unknown for
large $n$, and appears to be quite complicated. See \cite{bouc92a},
\cite{bouc92b} for low dimensional calculations.
Note for instance that the modules which occur do not seem to be limited to
any particular set of blocks of
$\F_p\Sigma_n$. It would be nice to have a better understanding of this
case. (See \cite{AF90} for a list of the
$p$-radical subgroups in
$\Sigma_n$.)
\enddemo

The next observation (which has also been used by group theorists \cite[p.\
479]{AS93}) gives, when combined with
Theorem~\ref{thm1}, a refinement of \cite[Prop.\ 6.1(iii)]{JMO92}.

\proclaim{Proposition} \label{coinvprop}
Let $K$ be a normal subgroup in $G$ of order prime to $p${\rm .}
 Then $|\pS|/K = |{\cal S}_p(G/K)|${\rm ,} and in particular $\St_*(G/K) =
\St_*(G)_K${\rm ,}
 where $(-)_K$ denotes
coinvariants{\rm .}\endproclaim

\demo{Proof}  To see the first claim we have to see that if $(P_0 \leq
\cdots \leq P_n)$ and $(P_0' \leq \cdots \leq P_n')$ are two chains of
subgroups in $G$ such that $P_iK = P_i'K$ for
all~$i$, then the two chains are conjugate by a single $k \in K$.  Since $K$
is normal in $G$ then by Sylow's theorem,
$PK = P'K$ if and only if we can find $k \in K$ such that $P' =kPk^{-1}$. In
our two chains above we can hence assume
that $P_n = P_n'$. But then $P_i$, $P_i'$ are subgroups of $P_n$ conjugate
by an element in $K$, so that $P_i = P_i'$
(since if $g
\in P_n, k
\in K$ and $kgk^{-1} \in P_n$ then $kgk^{-1}g^{-1} \in K \cap P_n = e$). The
second claim is a trivial consequence of
the first. 
\enddemo 

Although the homology groups of $\St_*(G)$ need not in general be projective
\cite{alperin90}, \cite{YW94} this is
however true for $H_0$, as shown by Quillen \cite[Prop.\ 5.2]{quillen78} and
Webb  \cite[pf.\ of Thm.\
E]{webb87comment}. Fix a Sylow $p$-subgroup $P$ of $G$ and let
$$\tilde G = \mbox{ subgroup of } G \mbox{ generated by }N_G(Q) \mbox{, }
Q\in \pS_{\leq P}.$$
(Using that $\pR$ and $\pS$ are $G$-equivalent and arguing as in
\cite[Prop.\ 5.2]{quillen78}, one sees that we may replace $\pS$ by $\pR$ in
the above definition
 of~$\tilde G$.)

\proclaimtitle{Quillen \cite{quillen78}, Webb \cite{webb87comment}}
\proclaim{Proposition} \label{h0formula}
$|\pS| = G \times_{\tilde G} |{\cal S}_p(\tilde G)|${\rm ,} and ${\cal
S}_p(\tilde G)$ is connected{\rm .} Furthermore
$H_0(\St_*(G)) = \ker(\Z_{(p)}[G/\tilde G] \pil{\rm aug} \Z_{(p)})$ is a
projective
$\Z_{(p)}G$\/{\rm -}\/module and $\St_*(G) \cong  H_0(\St_*(G)) \oplus
\St_*(\tilde G)\uparrow_{\tilde G}^G${\rm .}
\endproclaim 

\numbereddemo{{R}emark}\label{stronglypembedded}
Groups $G$ such that ${\cal S}_p(G)$ is disconnected are sometimes called
$p$-isolated \cite{goldschmidt70}. Groups such that ${\cal S}_p(G)$ is
disconnected or empty `have a strongly $p$-embedded subgroup'. (See e.g.,
\cite[Prop.\ 5.2]{quillen78}.) Such groups have been `classified'---see
\cite[(6.2)]{aschbacher93}.
\enddemo

\numbereddemo{Example}
If $G$ is a finite group of Lie type with $G \neq \tilde G$ then $G$ has
semisimple rank $1$. In this case there is, up to conjugacy, only one
nontrivial $p$-radical subgroup namely the Sylow $p$-subgroup and $\tilde G
= B$, the Borel subgroup.
The expression in Proposition \ref{h0formula} for $H_0(\St_*(G))$ is just
the standard alternating sum formula (cf.\ \cite[Thm.\ 66.35]{CR2}) for the
Steinberg module.
\enddemo

\proclaim{{C}orollary}  Let $\cC$ be a collection of $p$\/{\rm -}\/subgroups
closed under passage
to $p$\/{\rm -}\/overgroups and let $F: \bO_\cC^{\op} \to \Zpmod$ be a
functor concentrated  on conjugates of a subgroup
$P${\rm .}  Then 
$\lim^1_{\bO_\cC}F =
 \Hom_{WP}(H_0(\St_*(WP)),F(P))$ and $\lim^2_{\bO_\cC}F =
\Hom_{WP}(H_1(\St_*(WP)),F(P)).$
\endproclaim 

\demo{Proof} 
The statement for $\lim^1$ follows directly from Theorem~\ref{thm1} combined
with Proposition \ref{h0formula} above. The $\lim^2$ statement also follows
from these two results,
by the fact that $\Hom_G(-,M)$ is left exact.
\enddemo 

Note that this corollary has \cite[Prop.\ 6.2]{JMO92} as a consequence.

\numbereddemo{{R}emark} \label{lim1formula} If a projective module $P_S$
corresponding to a simple\break
$\F_pG$-module $S$ appears in $\SSt_0(G)$ then by Lemma
\ref{grplemma}, $S\downarrow^G_{\tilde G}$ contains a summand of the
projective $\F_p\tilde G$ module $P_{\F_p}$,
and $\dim_{\F_p}S \geq |G|_p$, a much stronger estimate than Proposition
\ref{stcrit}.

More precisely, using Frobenius reciprocity in the stable module category,
as in the proof of Lemma \ref{grplemma}, one sees that for any
$\F_pG$-module $M$
\begin{eqnarray*}
\dim_{\F_p}\Hom_G(\SSt_0(G;\F_p),M) &= &(\mbox{mult. of } P_{\F_p,\tilde G}
\mbox{ in }
M\downarrow^G_{\tilde G})\\[3pt]
&& -\  (\mbox{mult. of } P_{\F_p,G} \mbox{ in } M),
\end{eqnarray*} where $P_{\F_p,\tilde G}$ and
$P_{\F_p,G}$ denote the projective covers of $\F_p$ over $\F_p\tilde G$ and
$\F_pG$ respectively.
\enddemo

\numbereddemo{Example}\label{st0refinement}
Suppose that $g$ is an element of order $p$ in $G$ and that $P_S$ occurs in
$\SSt_0(G)$. Then $-\sum_{k=1}^{p-1}g^k$ does not act as the identity on
$S$, strengthening Proposition
\ref{stcrit}, Part $1$. To see this, note that we may assume $g \in \tilde
G$, so that $\F_p\langle g\rangle$ will be a
summand of $S \downarrow_{\langle g \rangle}$, and hence $1 +
\sum_{k=1}^{p-1}g^k = (g -1)^{p-1}$ does not act as
zero on $S\downarrow_{\langle g \rangle}$.
\enddemo

\numbereddemo{Example}Let $G = \Sigma_{2p}$ and suppose that $p \geq 5$.
Note that $\tilde G =  \Sigma_p \wr \Z/2$ and that $\tilde \Sigma_p = \Z/p
\rtimes \Z/(p-1)$. Let
$$M = \SSt_0(\Sigma_p) =\ker(\Z_p[\Sigma_{p}/(\Z/p \rtimes \Z/(p-1))]
\pil{\rm aug} \Z_p).$$
Note also that ${\cal S}_p(\tilde G) = {\cal S}_p(\Sigma_p) \star {\cal
S}_p(\Sigma_p)$, so that $\SSt_*(\tilde G) =
\SSt_1(\tilde G) = M\otimes_{\Z_p} M$, with the obvious $\Sigma_p \wr \Z/2$
action. Combination of these facts
describes the Steinberg complex of $\Sigma_{2p}$. Now,
$\SSt_0(G)  =\ker(\Z_p[\Sigma_{2p}/(\Sigma_p \wr \Z/2)] \pil{\rm aug}
\Z_p),$
$ \SSt_1(G)= (M \otimes_{\Z_p} M)\uparrow_{\Sigma_p \wr
\Z/2}^{\Sigma_{2p}},$
and the
differential is zero.
\enddemo

For later reference we also record the following useful fact, which follows
directly from the definition of $\SSt_*(G)$ as a suitably minimal summand of
a complex of permutation modules.

\proclaim{Proposition} \label{dual}
The complex $\SSt_*(G)$ is dimensionally self\/{\rm -}\/dual{\rm ,} i.e{\rm
., }
$\SSt_n(G) \cong \Hom_{\Z_p}(\SSt_n(G),\Z_p)$  as $\Z_pG$\/{\rm
-}\/modules{\rm .} 
\endproclaim  \pagebreak

\section{Higher limits over $\cC$-conjugacy categories} \label{conjugacysec}

The $\cC$-conjugacy category $\A_\cC$ is the category with objects
subgroups\break $H \in \cC$ and morphisms
homomorphisms $\varphi: H \to H'$ for which there exists $g \in G$ such that
$\varphi = g (\cdot)g^{-1}$. This
category occurs for example   when one works with the {\it
centralizer} homology decomposition (cf.\
\cite{JM92},
\cite{dwyer97}, \cite{dwyer98}). The
$\pSe$-conjugacy category is often called the {\it Frobenius category} of
$G$ and the $\pA$-conjugacy category is
often called the {\it Quillen category} of $G$.

Given a functor $F : \A_\cC \to \Zpmod$ we equip $|\cC|$ with the $G$-local
coefficient system $\cF$  given by
precomposition with the canonical map $(\dio \cC)_G \to \cC_G \to \A_\cC$
(see Remark \ref{grothrem}). On objects
we have $\cF(H_0 \leq \cdots \leq H_n) = F(H_n)$.

For higher limits of functors over $\A_\cC$ we obtain  results dual  to the
ones for $\bO_\cC$.

\proclaim{Theorem} \label{thm2} Let $\cC$ be any collection of nontrivial
$p$\/{\rm -}\/groups closed under
passage to
nontrivial elementary abelian subgroups{\rm ,} and let $F: \A_\cC \to
\Zpmod$ be any functor{\rm .} Then
$$\lim_{\A_\cC}{}^*F \cong H^*_G(|\cC|;\cF),$$
where $\cF$ is the $G$\/{\rm -}\/local coefficient system on $|\cC|$ induced
via
 $(\dio \cC)_G \to \cC_G \to \A_\cC${\rm .}

Furthermore for any collection $\cC'$ satisfying $\cC \cap \pA \subseteq
\cC' \subseteq \cC${\rm ,}
 the inclusion $|\cC'| \to |\cC|$ induces an isomorphism
$H^*_G(|\cC|;\cF) \pil{\simeq} H^*_G(|\cC'|;\cF)${\rm .}

If $F$ is a functor concentrated on a single conjugacy class of subgroups
with representative $V${\rm ,} then
$$
\lim_{ \A_\cC}{}^i F  =  \tuborg
\Hom_{NV/CV}(\St_{\GL(V)},F(V)) &\mbox{ if $V \in \pA$ and $i = \rk V - 1$}
\\ 0 & \mbox{otherwise}. \sluttuborg
$$
\endproclaim 

We prove this result later in this section.
(To see the duality to the case of orbit categories we note that
$\tilde C_*(|\pS_{< V}|)$ is chain homotopy equivalent to  $\St_{\GL(V)}$ in
degree $\rk V - 2$ if $V$ is elementary abelian and is contractible
otherwise.)
The formula for the higher limits for an atomic functor is due to Oliver
\cite{oliver94steinberg}, and using the rank filtration on $\pA$ we recover
his chain complex, arguing as in the proof of Theorem \ref{sscor}.

\proclaimtitle{\cite{oliver94steinberg}}
\proclaim{Theorem} \label{olivercomplex}\hskip-8pt
Under the assumptions and notation of Theorem~{\rm \ref{thm2}},
$\lim_{\A_\cC}^{*}F$ equals the homology of a
cochain complex $(C^*(F),\delta)$ with
$$C^i(F) =  \bigoplus_{[V] \in \cC \cap \pA/G; \rk V = i+1}
\Hom_{NV/CV}(\St_{\GL(V)},F(V)) .$$

The boundary map can be described as follows{\rm .}
Given $$\alpha \in \Hom_{NV/CV}(\St_{\GL(V)},F(V)),$$
 then the projection of $\delta(\alpha)$ onto
$\Hom_{NU/U}(\St_{\GL(U)},F(U))$  where $V \subset U$ and
$\dim_{\F_p}U = \dim_{\F_p}V +1$ is given by the composite

\begin{eqnarray*}
\noalign{\vskip-12pt}
 \St_{\GL(U)}& \to& \bigoplus_{V' \subset U; [V']=[V] \in \pA/G}
\St_{\GL(V')} \\[4pt]
&\pil{\alpha} &
 \bigoplus_{V' \subset U; [V']=[V]\in \pA/G} F(V') \to F(U).\end{eqnarray*}
 \endproclaim 

Before proving Theorem \ref{thm2} we establish some results.

\proclaim{Theorem}\label{thm2'}
Let $F: \A_\cC \to \Rmod$ be an atomic functor concentrated on conjugates of
$H \in \cC${\rm .} 
Then 
\begin{itemize} \ritem{1.}
$H^*_G(|E\A_\cC|;F) \cong H^*_{h(NH/CH)}(|\cC_{\leq
H}|,|\cC_{<H}|;F(H))$\hfill\break\hglue1in
$\cong \tilde H^*_{h(NH/CH)}(\Sigma |\cC_{<H}|;F(H)) .$
\ritem{2.}
If $R = \Z_{(p)}$ and $|\cC_{<H}|^Q$ is $\F_p$\/{\rm -}\/acyclic for all
nontrivial $p$\/{\rm -}\/subgroups in
 $NH/CH${\rm ,} then
 \begin{eqnarray*}
H^*_G(|E\A_\cC|;F)  \cong  H^*_G(|\cC|;\cF) &\hskip-7pt\cong &\hskip-7pt
H^*_{NH/CH}(|\cC_{\leq
H}|,|\cC_{<H}|;F(H))\\[4pt] &\hskip-7pt \cong&\hskip-7pt
\tilde H^*_{NH/CH}(\Sigma |\cC_{<H}|;F(H)) .
\end{eqnarray*}
\end{itemize}

\endproclaim 

\demo{Proof} 
The proof is dual to the proof of Theorem~\ref{thm1'}, and so we will be
brief. Writing out the definitions one sees
that:
$$ C^n_G(|E\A_\cC|;F) \cong C^n_{NH}(|{E\A_\cC}_{\leq
H}|,|{E\A_\cC}_{<H}|;F(H))$$
where ${E\A_\cC}_{\leq H}$ (resp.\ ${E\A_\cC}_{<H}$) denotes the full
subcategory of $E\A_\cC$ with objects the objects $f:H' \to G$ in $E\A_\cC$
satisfying $f(H')\leq H$ (resp.\ $f(H')< H$).
Hence we get corresponding identifications:
$$\begin{array}{ccc}
H^*_G(|\cC|;\cF) &\hskip.25in {\scriptstyle\cong}\hskip.25in  &
H^*_{NH}(|\cC_{\leq H}|,|\cC_{<H}|;F(H))\\[4pt]
\Big\downarrow&&\Big\downarrow\\[8pt]
H^*_G(|E\A_\cC|;F) & {\scriptstyle\cong} &  H^*_{NH}(|{E\A_\cC}_{\leq
H}|,|{E\A_\cC}_{<H}|;F(H)). \end{array}
$$
To show $1$ we note that $NH/CH$ acts freely on the pair $(|{E\A_\cC}_{\leq
H}|,|{E\A_\cC}_{<H}|)$ and argue as before. Likewise $2$ follows dually.
\enddemo 

The same proof as Corollary \ref{cor1'} now reveals:

\proclaim{{C}orollary}  \label{cor2'} Let $F: \A_\cC \to \Zpmod$ be a
functor and 
assume that $|\cC_{<H}|^Q$ is $\F_p$\/{\rm -}\/acyclic
for all $H \in \cC$ and all nontrivial $p$\/{\rm -}\/subgroups $Q$ of
$NH/CH${\rm .} Then
$H^*_G(|\cC|;F) \pil{\simeq} H^*_G(|E\A_\cC|;F)${\rm .}
\endproclaim 

{\it Proof of Theorem~{\rm \ref{thm2}}}.
We already know that $\lim^*_{\A_\cC} F = H^*_G(|E\A_\cC|;F)$ from
Propositions \ref{prop1} and \ref{prop2}, and so
to prove the first part of the theorem, we just have to show that the
conditions of Corollary \ref{cor2'} apply.

We start by assuming that $\cC$ is closed under passage to all nontrivial
subgroups, and will later show that nonelementary abelian $p$-subgroups can
be removed without changing either higher limits or Bredon cohomology.
Hence, for all $H \in \cC$, $\cC_{<P} = { \cal S}_p(G)_{<P}$.
If $Q$ is a nontrivial $p$-subgroup of $NP/CP$, then $|\pS_{<P}^Q|$ is
contractible by the inequalities $H \geq H^Q \leq P^Q$. (The
fixed-point subgroup $H^Q$ is not trivial since $Q$ and $H$ are $p$-groups.)
Now Corollary~\ref{cor2'} and
Theorem~\ref{thm2'} prove Theorem~\ref{thm2}, except the statement about
$\cC'$ and the last rewriting in the
atomic case.

We address the last rewriting first. If $P$ is not elementary abelian, then
the inequalities $H \leq H\Phi(P) \geq
\Phi(P)$ show that $|{\cal S}_p(G)_{<P}|$ is $NP$-equivariantly contractible
(see also \cite{TW91}). (Here $\Phi(P)$
denotes the Frattini subgroup of $P$; alternatively, we use the inequalities
$H \geq H \cap [P,P] \leq [P,P]$ if $P$ is
non-abelian, and the inequalities $H \geq {}_pH \leq {}_pP$  if $P$ is
abelian, instead.) Hence $\tilde C_*(|\pS_{<P}|)$
is $NP$-contractible if $P$ is not elementary abelian. If $P=V$ is
elementary abelian then $|\pS_{<V}|$ is the space
associated to the ordered simplicial complex of flags of nontrivial proper
sub-vector spaces of $V$, which again is
equal to the barycentric subdivision of the Tits building of $\GL(V)$, by
the definition of the Tits building
\cite{tits74}. But $\tilde C_*(|\pS_{<V}|)$ is $NV$-chain homotopy
equivalent to its homology, which is just the
Steinberg representation of $\GL(V)$ in dimension $\rk V - 2$. (One way to
see this is to use Webb's
Theorem~\ref{webbthm}.) This shows the last rewriting in the atomic case.

To see that we can replace $\cC$ by $\cC'$ we use the same argument as in
the proof of Theorem~\ref{thm1}. We can suppose that $F$ is an atomic
functor concentrated on conjugates of a subgroup $P \in \cC$. If $P \in \cC
\setminus \cC'$ then we are done by the formula for an atomic functor just
proved. If $P \in \cC'$ then we are done if we can show that $\cC_{<P}$ and
$\cC'_{<P}$ are $NP$-homotopy equivalent. This is done by showing that we
can remove subgroups in $\cC_{<P} \setminus \cC'_{<P}$ from $\cC$ in order
of decreasing size, without changing the $NP$ homotopy type, using the
pushout technique of the proof of Theorem~\ref{thm1} together with the
$NP'$-contractibility of $\pS_{<P'}$ when $P'$ is not elementary abelian, as
shown above.
\hfill\qed\vglue4pt
\pagegoal=49pc

{\it Remark} 6.5.
A natural question to consider is whether one can obtain a similar
description of higher limits of functors over
quotient categories~$\D$ of $\cC_G$ or $(\cC^{\op})_G$ other than just
$\A_\cC$ and $\bO_\cC^{\op}$. The proofs of
Theorems~\ref{thm1'} and \ref{thm2'} show what has to be satisfied. For
quotient categories of $\cC_G$ one needs
that $|\cC_{<H}|^Q$ is $\F_p$ acyclic for all $H \in \cC$ and all nontrivial
$p$-subgroups $Q \leq \Aut_\D(H)$ (and
analogously for $(\cC^{\op})_G$). We see that this property is rather
special for $\A_\cC$ (and $\bO_\cC^{\op}$), and
generally does not hold for, for example, their opposite categories. For
their opposite categories we however get
formulas like the ones found in this paper for higher colimits.

\section{Higher limits over orbit simplex categories} \label{simplexsec}
\pagegoal=48pc

Let $(\sd\cC)/G$ be the orbit simplex category, i.e., the poset with objects
$G$-conjugacy classes {$[\sigma]$ of
proper chains of subgroups $\sigma = (H_0 < \cdots < H_n)$ and with a
morphism $[\sigma] \to [\sigma']$ if we can
find representatives $\sigma,\sigma'$, such that $\sigma$ is a refinement of
$\sigma'$. This category occurs when
we work with the {\it normalizer} homology decomposition (cf.\
\cite{dwyer97}, \cite{dwyer98}). Higher limits of
functors on $((\sd\cC)/G)^{\op}$ can likewise be expressed as equivariant
cohomology of the subgroup complex $|\cC|$.
Indeed, this follows by a formal argument, and does not require any
assumptions on the collection $\cC$ or the ground
ring.

Since $(\sd \cC)/G$ is already a finite poset and hence  $|(\sd \cC)/G| =
|\sd \cC|/G$ is a finite space (with associated topological space
homeomorphic to $|\cC|/G$), the only improvement will be to remove the
subdivision $\sd$. This follows easily from an argument of  S\l omi\'nska
\cite{slominska92}. (Alternatively one can use Remark \ref{sheafremark}: The
topological spaces associated to $|\sd \cC|$ and $|\cC|$ are
$G$-homeomorphic, and under the canonical homeomorphism the sheaves on the
topological spaces associated to $|\sd \cC|/G$ and $|\cC|/G$ coincide.)
The following proposition is due to her.

\proclaim{Proposition} \label{normprop}\hskip-8pt
 Let $\cC$ be a  collection of subgroups of a finite group~$G${\rm ,} and
let $F:((\sd\cC)/G)^{\op}
\to \Rmod$ be an arbitrary functor{\rm .} Then
$$\lim_{(\sd \cC)/G}\hspace{-10pt}{}^*F \cong H^*_G(|\cC|;\cF) ,$$
where $\cF$ is the $G$-local coefficient system given via $(\dio \cC)_G \to
((\sd \cC)/G)^{\op}${\rm .}
\endproclaim 

\demo{Proof}  Note that $|\cC|$ is in fact an ordered simplicial complex, so
that the subcomplex $|\sigma|$ generated by an $n$-simplex $\sigma$ is
contractible. Hence the functor $\sigma \mapsto C_*(|\sigma|;R)$, i.e., the
functor which to $\sigma$ associates the normalized chain complex of the
$n$-simplex $\sigma$, gives a more economical projective resolution of $R$
in $\Rm{(\sd \cC)/G}$, than the one of Proposition \ref{prop1}. In this case
we see that $\Hom_{R{(\sd \cC)/G}}(C_n(|-|;R),F) = \prod_{[\sigma]}
F(\sigma)$, where $[\sigma]$ runs over the $G$-conjugacy classes of
nondegenerate $n$-simplices in $|\cC|$. Hence $\lim^*_{(\sd \cC)/G}F =
H^*(|\cC|/G;F) = H^*_G(|\cC|;\cF)$ as wanted.
\enddemo   

In the topological applications occurring when one uses the normalizer
decomposition, $F$ will most often be induced
by a generic coefficient system $\fF :{\bO}(G)^{\op} \to \Rmod$; i.e., $F$
is the functor which on objects is given by
$[(P_0 \leq \cdots \leq P_n)] \mapsto \fF(NP_0 \cap \cdots \cap NP_n)$. In
this case we have $C^*_G(|\cC|;\cF) =
C^*_G(|\cC|;\fF)$. (Note however that $\cF$ is not exactly the $G$-local
coefficient system induced by $\fF$.)  Hence
in this case we may replace $|\cC|$ with any $G$-homotopy equivalent space,
for the purpose of calculating
$H^*_G(|\cC|;\cF)$. More precisely we get:

\proclaim{{C}orollary}  \label{easycor} Let $\cC$ be a collection of
subgroups of a finite group $G$ and let $F:
(\sd\cC)/G \to \Rmod$ be a functor induced by a generic coefficient system
$\fF${\rm .} Then $$\lim_{(\sd
\cC)/G}\hspace{-10pt}{}^*F \cong H^*_G(|\cC'|;\fF)$$  for any collection
$\cC'$ such that $|\cC'|$ is $G$\/{\rm -}\/homotopy
equivalent to $|\cC|${\rm .}
\endproclaim 

We end this section by examining the special case where $\fF$ is the generic
coefficient system $\fF(-) =
H^n(-;\Z_{(p)})$, $n\geq 0$, which illustrates that for a  particular $\fF$
one can do better than
Corollary~\ref{easycor}.

\proclaim{Theorem}\label{normsharp}
Let $\cC$ be a collection of $p$\/{\rm -}\/subgroups of a finite 
group~$G${\rm ,}
 closed under passage to $p$\/{\rm -}\/overgroups{\rm .}
 Let $\fF(-) = H^n(-;\Z_{(p)})$, $n \geq 0${\rm ,} and let 
$F: ((\sd \cC)/G)^{\op} \to \Zpmod$ be the
associated
functor on objects given by $[(P_0 \leq \cdots \leq P_n)] \mapsto \fF(NP_0
\cap \cdots \cap NP_n)${\rm .}
 Then
$$  \lim_{(\sd \cC)/G}\hspace{-10pt}{}^*F \cong \lim_{(\sd
\cC')/G}\hspace{-10pt}{}^*F \cong H^*_G(|\cC'|;\fF),$$
 for $\cC'$ any collection $\cC \cap \pD \subseteq \cC' \subseteq \cC${\rm
.}
\endproclaim 

\demo{Proof}  As in the proof of Theorem~\ref{pcentriccor} it is enough to
prove the statement over $\Z_p$, and we therefore adopt the convention that
all coefficients are
over~$\Z_p$. We want to show that we can
successively remove conjugacy classes of subgroups in $\cC \setminus \cC'$
in order of increasing size, without
changing the Bredon cohomology. So suppose that $ \tilde \cC$ is a
subcollection obtained from $\cC$ by removing
$G$-conjugacy of subgroups in this way and let $\hat \cC$ denote the
collection obtained from $\tilde \cC$ by
removing the $G$-conjugacy class of a minimal subgroup $P \in \tilde \cC
\setminus \cC'$. As in the proof of
Theorem~\ref{thm1} we have a pushout square, which is a homotopy pushout of
$G$-spaces 
$$
\begin{array}{ccc} G \times_{NP}|\link_{\tilde \cC}P| &\lrar & |\hat \cC|
\\[4pt]
       \Big\downarrow&&\Big\downarrow\\[8pt]
    G \times_{NP}|\sta_{\tilde\cC}P|&\lrar & |\tilde \cC| .\end{array} $$
We need to see that $\tilde H^*_{NP}(|\link_{\tilde \cC}P|;\fF) = 0$.
Writing out the definition, we must see that $H^n(NP;\tilde
C_*(|\link_{\tilde \cC} P|))$ is an acyclic chain complex.
But as in the proof of Theorem~\ref{pcentriccor}, $\tilde C_*(|\link_{\tilde
\cC}P|)$ is an $NP$-chain homotopy
equivalent to 
$\SSt_*(WP) \otimes \Sigma\tilde C_*(|\tilde \cC_{<P'}|)$.

Consider the Lyndon-Hochschild-Serre spectral sequence:
\begin{eqnarray*}
 E_2^{i,j} &= &H^i(WP;H^j(P;\tilde C_r(|\tilde \cC_{<P}|))
\otimes \SSt_m(WP) ) \\
&\Rightarrow& H^{i+j}(NP; \SSt_m(WP) \otimes \tilde C_r(|\tilde \cC_{<P}|)).
\end{eqnarray*}
Since $H^j(P;\tilde C_r(|\tilde \cC_{<P}|))\otimes \SSt_m(WP)$ is a
projective $WP$-module, the spectral sequence collapses and yields
the fact that the chain complex
$$H^n(NP;\tilde C_*(|\link_{\tilde \cC}P|))$$ is chain homotopy equivalent
to the complex
$$\Hom_{WP}(\Z_p, H^n(P;\tilde C_*(|\tilde \cC_{<P}|))\otimes \SSt_*(WP)).$$
But, using Proposition \ref{dual}, we
obtain
\begin{eqnarray*}
&&\hskip-.5in\Hom_{WP}(\Z_p, H^n(P;\tilde C_r(|\tilde \cC_{<P}|))\otimes
\SSt_m(WP)) \\[4pt]
&&\qquad\quad\cong \Hom_{WP}(
\SSt_m(WP),H^n(P;\tilde C_r(|\tilde \cC_{<P}|)))
\end{eqnarray*}
 which is zero if $P$ is non-$p$-centric by
Corollary~\ref{stcritcor}. It is also zero if $P$ is\break
$p$-centric, since in this case it equals  $$\Hom_{WP}(
\SSt_m(NP/PC(P)),H^n(P;\tilde C_r(|\tilde \cC_{<P}|))),$$ by
Proposition \ref{coinvprop}, which is zero since $O_p(NP/PC(P)) \neq e$.
Hence in both cases our complex
$H^n(NP;\tilde C_*(|\link_{\tilde \cC} P|))$ is chain homotopy
 equivalent to the trivial complex, and therefore in particular acyclic,
which finishes the proof of the theorem.
\enddemo 

\numbereddemo{{R}emark}
The key feature in the above argument was being able to replace
$\St_*(WP;\Z_p)$ with the chain homotopic complex $\SSt_*(WP)$. This will
be possible  more generally  if just $F$ extends from a functor
${\bO}(G)^{\op} \to \Zpmod$ to a $\Z_{(p)}$-linear
functor $\Zpm{G}
\to \Zpmod$. (Having an extension to a $\Z_{(p)}$-linear functor
$\Z_{(p)}G$-(permutation modules) $\to \Zpmod$ is
equivalent to $F$ being a cohomological Mackey functor, by
\cite{yoshida83}.) Note also that in concrete situations the
method of the above proof can be used to determine if further subgroups can
be removed.
\enddemo

{\it Example} 7.5.
To give an example where one cannot get away with just considering $\pD$,
consider the generic coefficient system $F(H) = R(H) \otimes \Z_{(2)}$,
where $R$ 
denotes the complex representation ring. (This is a non-cohomological Mackey
functor.) If $G = \Sigma_3 \times \Z/2$
then $H^*_G(|{\cal S}_2(G)|;F) = R(G) \otimes \Z_{(2)}$ and  $H^*_G(|{\cal
D}_2(G)|;F) = R(\Z/2 \times \Z/2) \otimes
\Z_{(2)}$. These are obviously different, the first having rank $6$, the
latter $4$. (See \cite{thevenaz93} for
interesting work using this coefficient system.)
 
\section{Induction sequences} \label{exactseqsec}

The purpose of this section is to examine our results in the special case
where $F$ is a functor for which all the higher limits in fact vanish.
In this case our results produce exact sequences expressing $F(G)$ in terms
of $F$ evaluated on subgroups,  the one corresponding to the normalizer
decomposition being the one found by Webb \cite{webb91}.
We do not assume that our functors are Mackey functors, although a Mackey
structure is often needed in order to show
that the higher limits indeed vanish.

We expand on work of Dwyer \cite{dwyer98} and Webb \cite{webb91} (see also
Dress \cite{dress75}), the main new results being
 Examples~8.8 and~8.9 and Corollaries \ref{exactthm1}, \ref{exactthm2},
and~\ref{exactliecor}.

Our general setup is the following. Suppose we have a $G$-space $X$, a
functor $F: ({\cal S}(G)^{\op})_G \to \Zmod$, and a map $(\dio X)_G \to
({\cal S}(G)^{\op})_G$ making $F$ a $G$-local coefficient system on $X$. If
$H^n_G(X;F) = 0$ for $n>0$ and $H^0_G(X;F) = F(G)$, then writing
the normalized cochain complex gives us an exact sequence
$$ 0 \to F(G) \to \prod_{\sigma_0 \in X_0/G} F(\sigma_0)^{G_{\sigma_0}} \to
\prod_{\sigma_1 \in X_1/G} F(\sigma_1)^{G_{\sigma_1}} \to \cdots ,$$
by setting $C^{-1}(X;F) = F(G)$.
Here $X_i$ denotes the set of nondegenerate $i$-simplices of $X$.
This sequence expresses $F(G)$ in terms of subgroups in a quite explicit
way. We will describe some instances of the above situation, as well as some
related methods.

\demo{{\rm 8.1.} Webb\/{\rm '}\/s exact sequence} We will first try to
explain the work of Webb in the current context---we
refer the reader to \cite{webb91}, \cite{webb99} for much more information,
as well as many beautiful applications.

We say that a functor $\fF: {\bO}(G)^{\op} \to \Rmod$ is projective (as a
generic coefficient system) relative to a collection $\cC$ if the canonical
map $\oplus_{H \in \cC}\fF\downarrow^G_H\uparrow_H^G \to \fF$, 
given as the unit of an adjunction, is split surjective in the
category of functors  ${\bO}(G)^{\op} \to \Rmod$. (For example if $\cC$ is a
cohomological Mackey functor to $\Zpmod$ then a transfer argument shows that
$\fF$ is projective relative to $\cC = \pS$.)
\advance\theoremcount by 1

\proclaimtitle{Webb \cite{webb91}}
\proclaim{Theorem} \label{webbseqthm}
 Let $\fF: {\bO}(G)^{\op} \to \Rmod$ be any functor{\rm .} Let $X$ be a
$G$\/{\rm -}\/space{\rm .}
Consider two collections  $\cC' \subseteq \cC$ of subgroups
of $G$ both closed under passage to subgroups{\rm .} Assume that $\fF$ is
projective
 relative to $\cC${\rm ,} that $\fF(H) = 0$
for all $H \in \cC'${\rm ,} and that $X^H$ is contractible for all $H \in
\cC \setminus \cC'${\rm .}

Then $H^*_G(X;\fF) = \fF(G)$ and there is a split exact sequence
$$ 0 \to \fF(G) \to \prod_{\sigma_0 \in X_0/G} \fF(G_{\sigma_0})
 \to \prod_{\sigma_1 \in X_1/G} \fF(G_{\sigma_1}) \to \cdots.$$
\endproclaim 

\demo{Proof}  Since $\oplus_{H \in \cC}\fF\downarrow^G_H\uparrow_H^G \to
\fF$ is assumed to be split surjective it is enough to prove the statements
for $\fF\downarrow^G_H\uparrow_H^G$. However, since
for any $H$-generic coefficient system $\fF$ we have the adjunction
$C^*_G(X;\fF\uparrow_H^G) = C^*_H(X;\fF)$ we can without restriction assume
that $\cC$ is the collection of all subgroups. The general acyclicity result
now follows from Lemma \ref{genericlemma} below,
in the special case when $Y = *$. That the sequence in the theorem is in
fact split exact follows from the proof of
Lemma \ref{genericlemma}, since we may replace $X$ by the $G$-contractible
subspace $X'$.
\enddemo 
 
\proclaim{Lemma} \label{genericlemma} Let $f: X \to Y$ be a map of $G$\/{\rm
-}\/spaces{\rm ,}
 $\cC$ a collection of subgroups of $G$ closed under passage to
subgroups{\rm ,} and $\fF:{\bO}(G)^{\op} \to \Rmod$ an
arbitrary generic coefficient system{\rm .} Assume that $F(H) = 0$ for all
$H \in \cC$ and that $f^H: X^H \to Y^H$ is a
homotopy equivalence for all $H \not \in \cC${\rm .} Then $H^*_G(Y;\fF)
\pil{\simeq} H^*_G(X;\fF)${\rm .}
\endproclaim

\demo{Proof}  Define $X' = \cup_{K\not \in \cC}X^K$ and note that
$C^*_G(X;\fF) = C^*_G(X';\fF)$ by the assumption on
$\fF$. We define $Y'$ similarly. It hence suffices to show that the
restriction of $f$ to a map $f': X' \to Y'$ is a $G$-equivalence. If $H
\not \in \cC$ then ${X'}^H = X^H$ and ${Y'}^H =Y^H$ since $\cC$ is closed
under passage to subgroups, so ${f'}^H$ is a
homotopy equivalence by our assumption. If $H \in \cC$ we need to see that
${f'}^H: {X'}^H = \cup_{H\leq K, K \not \in
\cC}X^K \to \cup_{H\leq K, K \not \in \cC}Y^K ={Y'}^H$ is a homotopy
equivalence. But this is true since on each of the
pieces $X^K$ and on their intersections which are of the same form, the map
${f'}^H_{|X^K}: X^K \to Y^K$ is a homotopy
equivalence. (The most appealing argument for why this implies that ${f'}^H$
is a homotopy equivalence is perhaps to
interpret the union as a homotopy colimit, and refer to the homotopy
invariance property of homotopy colimits.) Hence
$f'$ is a homotopy equivalence on all fixed points, and so we conclude that
$f'$ is a $G$-equivalence.
\enddemo 

It is interesting to compare the method of the above lemma to the `method of
discarded orbits' in \cite[\S 6.7]{dwyer98}.

\demo{{\rm 8.4.} Exact sequences arising from an acyclic collection}
Let $\cC$ be a collection of subgroups in $G$ and let $\fF: {\bO}(G)^{\op}
\to \Zpmod$ be any functor. To $\fF$ we can associate three other functors.
Let $F^{\beta}: \bO_\cC^{\op} \to \Zpmod$ be the functor given by $H \mapsto
\fF(H)$. Let $F^{\alpha}: \A_\cC \to \Zpmod$ be the functor given by $H
\mapsto \fF(CH)$. Finally let $F^{\delta}: ((\sd \cC)/G)^{\op} \to \Zpmod$
be the functor given by $[H_0 < \cdots < H_k] \mapsto \fF(NH_0 \cap \cdots
\cap NH_k)$.
Note that we have canonical maps $\fF(G) \to \lim^0_{\bO_\cC}F^\beta$, $F(G)
\to \lim^0_{\A_\cC}F^\alpha$, and  $\fF(G) \to \lim^0_{(\sd \cC)/G}F^\delta$
given by restriction.

\advance\theoremcount by 1
\numbereddemo{Definition} We say that $\fF$ is {\it subgroup $\cC$-acyclic}
if $\lim^i_{\bO_\cC} F^\beta = 0$ for $i>0$ and $\fF(G) \pil{\simeq}
\lim^0_{\bO_\cC}F^\beta$. Likewise we say that $\fF$ is {\it centralizer or
normalizer $\cC$-acyclic} if the analogous conditions hold for $F^\alpha$ or
$F^\delta$.
\enddemo

By the rewritings of Propositions \ref{prop1}, \ref{prop2}, and
\ref{normprop}, this definition can be stated as stating
that $\fF$ is subgroup (resp.\ centralizer and normalizer) $\cC$-acyclic if
and only if the $G$-space $|E\bO_\cC|$
(resp.\
$|E\A_\cC|$ and $|\cC|$) is acyclic with respect to the generic coefficient
system $\fF$.

We now give some important examples of acyclicity before giving our exact
sequences.

\numbereddemo{Example}[The collection $\pSe$]\label{collpSe}
Consider the case $\cC = \pSe$, and let $M$ be an arbitrary
$\Z_{(p)}G$-module.
By the analysis of Section~\ref{prelim}, the spaces  $|E\bO_\cC|$,
$|E\A_\cC|$ and $|\cC|$ are $P$-contractible, where $P$ is a Sylow
$p$-subgroup of $G$.
Since $H^n(-;M)$, $n\geq 0$, is a cohomological Mackey functor, a transfer
argument shows that $H^n(-;M)$ is both subgroup,
centralizer  and normalizer $\pSe$-acyclic.
\enddemo

{\it Example} 8.6 (The collection $\pS$). Let $\fF(-) = H^n(-;M):
{\bO}(G)^{\op} \to \Zpmod$, where $M$ is some
$\Z_{(p)}G$-module.
 By our analysis in Sections \ref{orbitsec} and~\ref{conjugacysec}, the
spaces  $|E\bO_{\pS}|$,
$|E\A_{\pS}|$ and $|\pS|$
are $P$-equivalent. Hence by an application of the transfer, we see that
$$H^*_G(|E\bO_{\pS}|;\fF) \cong
H^*_G(|{\pS}|;\fF) \cong  H^*_G(|E\A_{\pS}|;\fF) .$$ For $n>0$ the formula
in Theorem~\ref{thm1} shows that
$$H^*_G(|E\bO_{\pS}|;\fF) = H^*_G(|E\bO_{\pSe}|;\fF),$$ since $H^n(e;M) = 0$
for $n>0$. So, by the previous example, we
see that  $H^n(-;M)$, $n > 0$, is both subgroup, centralizer and normalizer
$\pS$-acyclic. (Alternatively, a proof can be
obtained by appealing to Lemma \ref{genericlemma}.) For $n=0$, by
definition,   $H^*_G(|\pS|;M) = M^G$ if and only if
$\Hom_G(\St_*(G),M)$ is an acyclic cochain complex.   Hence $H^n(-;M)$ is
subgroup $\pS$-acyclic if and only if it is
centralizer
$\pS$-acyclic if and only if it is normalizer $\pS$-acyclic, iff
$\Hom_G(\St_*(G),M)$ is an acyclic cochain complex.
\vglue12pt

By removing subgroups we can now establish sharpness of smaller collections.

\vglue12pt {\it Example} 8.8 (The collection $\pD$). Example \ref{collpSe}
combined with
Theorem~\ref{pcentriccor} and
\ref{normsharp} and Corollary \ref{easycor} shows that  $\fF =
H^n(-;\Z_{(p)})$ is subgroup and normalizer
$\pD$-acyclic, $n\geq 0$.

\vglue12pt {\it Example} 8.9 (The collections $\pDeM$ and $\pDM$).  The
previous example can be
generalized to deal with arbitrary coefficient modules. Let $M$ be a
$\Z_{(p)}G$-module and let $K$ be the kernel of
the homomorphism $G \to \Aut M$. As in \cite{dwyer97} a $p$-subgroup $P$ is
said to be $M$-centric if $Z(P) \cap K$
is a Sylow
$p$-subgroup in $C(P) \cap K$. (If $P$ is $M$-centric then we can write
$C(P) \cap K = (Z(P) \cap K) \times L$ for
some subgroup $L$ where $p \hbox{$\not|$}\,  L|$.) Let $\pDeM$ denote the
collection of $M$-centric subgroups which satisfy
$O_p(NP/(P(C(P)\cap K))) =e$, and let $\pDM = \pDeM \cap \pS$. Note that
$\pD \subseteq \pDeM \subseteq \pRe$. In
\pagebreak
the extreme case where  $M$ is the trivial module and $p | |G|$ then $\pDM =
\pDeM = \pD$, since $e$ is noncentric. On
the other hand if $G$ acts faithfully on $M$, then $\pDeM = \pRe$.

The proofs of Theorems~\ref{pcentriccor} and \ref{normsharp} reveal that
$H^n(-;M)$ is subgroup and normalizer
$\pDeM$-acyclic for $n\geq 0$. Likewise $H^n(-;M)$ is subgroup and
normalizer $\pDM$-acyclic for $n>0$, and for $n=0$ it is $\pDM$-acyclic if 
and only if it is $\pS$-acyclic.
 
\advance\theoremcount by 3

\numbereddemo{{R}emark}
By a result of Diethelm-Stammbach \cite{DS84} (see also \cite{HK88},
\cite{symonds99}), for any finite $p$-group
$P$ and any $W \leq \Out(P)$, every simple\break $\F_pW$-module appears as a
composition factor of $H^n(P;\F_p)$
for some
$n$.  Hence properties of $\St_*(W)$ will not in the most naive way allow
one to remove more subgroups from the
collection $\pD$ when examining acyclicity properties of the functors
$H^n(-;\Z_{(p)})$, if one wants the collection to
work for {\it all} $n \geq 0$. (But, see Example \ref{miyamoto}.) It would
be interesting to get smaller collections
than $\pDeM$ for the functors $H^n(-;M)$, using more elaborate information
about $M$. 
\enddemo 

\numbereddemo{{R}emark}
For the subgroup and normalizer $\pDM$-acyclicity of $H^n(-;M)$, $n>0$ we
can furthermore remove those subgroups $P \in \pDeM$ such that
$M\downarrow^G_P$ is projective. This is obvious for subgroup acyclicity,
and for normalizer acyclicity it follows from Lemma \ref{genericlemma}.
\enddemo

\numbereddemo{Example}[The collection $\pA$]\label{pAexample}
Theorem~\ref{thm2} and Corollary \ref{easycor} show that $H^n(-;M)$ is
centralizer and normalizer $\pA$-acyclic for $n>0$ and has the same
centralizer and normalizer $\pA$-acyclicity properties as 
$\pS$-acyclicity properties for $n=0$.
\enddemo 

Theorem~\ref{thm1}, Theorem~\ref{thm2}, and Proposition \ref{normprop}
immediately produce three exact sequences, where the last one is a version
of Webb's result.

\proclaim{{C}orollary}  \label{exactthm1} Let $\cC$ be a collection of
$p$\/{\rm -}\/subgroups of a finite group $G$,
 closed under passage to $p$\/{\rm -}\/radical overgroups{\rm .}
 Let $\fF: {\bO}(G)^{\op} \to \Zpmod$ be any subgroup $\cC$\/{\rm
-}\/acyclic
functor{\rm .} Then there is an exact sequence
$$ 0 \to \fF(G) \to \oplus_{[P] \in \cC'/G} \fF(P)^{NP} \to \oplus_{[P_0 <
P_1] \in |\cC'|_1/G} \fF(P_0)^{NP_0
\cap NP_1} \to \cdots \to 0$$
where $\cC' = \cC \cap \pRe${\rm .}
 The first boundary map is given by restriction and the latter comes from
the boundary map in the cochain complex
$C^*_G(|\cC'|;F^{\beta})${\rm .}

If for all $P \in \cC${\rm ,} $CP$ acts trivially on $\fF(P)$ then  $\cC'$
may be replaced by $\cC \cap \pD${\rm .}
\endproclaim 

\proclaim{{C}orollary}  \label{exactthm2} Let $\cC$ be a collection of
nontrivial $p$\/{\rm -}\/subgroups
 of a finite group $G$, closed under passage to nontrivial elementary abelian
subgroups{\rm .}
  Let $\fF: {\bO}(G)^{\op} \to \Zpmod$ be any
centralizer $\cC$\/{\rm -}\/acyclic functor{\rm .} Then there is an exact
sequence
$$ 0 \to \fF(G) \to \oplus_{[V] \in \cC'/G} \fF(CV)^{NV} \to \oplus_{[V_0 <
V_1]
 \in |\cC'|_1/G} \fF(CV_1)^{NV_0\cap NV_1} \to \cdots \to 0,$$
where $\cC' = \cC \cap \pA${\rm .}
The first boundary map is given by restriction
and the latter comes from the boundary map in the cochain complex
$C^*_G(|\cC'|;F^{\alpha})${\rm .}
\endproclaim 

\proclaimtitle{\cite{webb91}}
\proclaim{{C}orollary} \label{exactthm3}Let $\cC$ be any collection of 
subgroups of a finite group $G${\rm .}
Let  $\fF:
{\bO}(G)^{\op}
\to
\Rmod$ be any normalizer $\cC$\/{\rm -}\/acyclic functor{\rm .} Then there
is an exact sequence
$$ 0 \to \fF(G) \to \oplus_{[H]
 \in \cC/G} \fF(NH) \to \oplus_{[H_0<H_1] \in |\cC|_1/G} \fF(NH_0\cap NH_1)
\to \cdots \to 0,$$
where the boundary map comes from the boundary map in the cochain complex
$\tilde C^*_G(|\cC|;\fF)$, where $\fF$ is viewed as a generic coefficient
system{\rm .}
\endproclaim 

\numbereddemo{Example} Consider the sporadic group $M_{11}$. According to
the Atlas \cite{atlas} there are two conjugacy classes of elementary abelian
$2$-subgroups:  one of rank one $V_1$ with centralizer $CV_1 = \GL_2(\F_3)$
and one of rank two $V_2$ which is
self-centralizing $CV_2=V_2$. Furthermore $NV_1 = CV_1$, $NV_2 = \Sigma_4$,
and $NV_1 \cap NV_2 = D_8$.
Example \ref{pAexample} and Corollary \ref{exactthm2} hence produce an exact
sequence for any $n>0$ and any
$\Z_{(2)}M_{11}$-module $N$,
$$0 \to H^n(M_{11};N) \to H^n(\GL_2(\F_3);N) \oplus H^n(V_2;N)^{\Sigma_3}
\to H^n(V_2;N)^{\Z/2} \to 0 .$$
\enddemo

The chain complexes of Theorem \ref{liecor} and Theorem \ref{olivercomplex}
likewise give  us exact sequences.

\proclaim{{C}orollary}  \label{exactliecor} \hskip-8pt
Let $G$ be a finite group of Lie type of characteristic~$p${\rm ,}
 and let $\cC$ be a collection of $p$\/{\rm -}\/subgroups closed under
passage to $p$\/{\rm -}\/radical overgroups{\rm .}
 Let $\fF:
{\bO}(G)^{\op} \to \Zpmod$ be any subgroup $\cC$\/{\rm -}\/acyclic
functor{\rm .} Then there is an exact sequence
$$0 \to \fF(G) \to  \fF(U_{\emptyset})^{P_\emptyset} \to \bigoplus_{J
\subseteq S ; |J|=1
 ; U_J \in \cC} \Hom_{L_J}(\St_{L_J},\fF(U_J)) \to \cdots \to 0.$$
The first boundary map is restriction to the Sylow $p$\/{\rm -}\/subgroup
 $U_{\emptyset}$ in $G$  and the following ones are described in
Theorem~{\rm \ref{identprop}.}
\endproclaim 

\proclaim{{C}orollary}
Let $\cC$ be a collection of nontrivial $p$\/{\rm -}\/subgroups
of a finite group $G$, closed under passage to nontrivial elementary abelian
$p$\/{\rm -}\/subgroups{\rm .}
 Let $\fF: {\bO}(G)^{\op} \to \Zpmod$ be any centralizer $\cC$\/{\rm
-}\/acyclic functor{\rm .}
 Then we have an exact
sequence
\begin{eqnarray*}
0 &\to &\fF(G) \to   \bigoplus_{[V] \in \cC \cap \pA/G; \rk V = 1}
\hspace{-20pt}\fF(CV)^{NV}\\[6pt] & \to   &\bigoplus_{[V] \in \cC \cap
\pA/G; \rk V = 2} 
 \Hom_{NV/CV}(\St_{\GL(V)},\fF(CV)) \to \cdots \to 0.
\end{eqnarray*}

The first boundary map is given by restriction and the following ones are
described in Theorem {\rm \ref{olivercomplex}.}
\endproclaim 

Theorem \ref{sscor} can also be used to obtain information about $\fF(G)$.
However, in general, if $G$ is not a finite group of Lie type, then we do
not get an exact sequence, but we can get two other expressions. One is
obtained by looking at what happens at the $(0,0)$-spot in the spectral
sequence, where we get a formula generalizing
the Cartan-Eilenberg formula for stable elements in group cohomology
\cite{CE56} to be   explained   in
Section~\ref{stableelements}. Another is obtained by trading in the spectral
sequence for an expression in the
Grothendieck group of modules, where we get a somewhat generalized version
of a result of Webb, which we now
explain.

Define the Lefschetz module $L$ of a chain complex $C_*$ to be the virtual
module $L(C_*) = \sum_i(-1)^iC_i$, where the alternating sum is taken in the
Grothendieck group of modules under direct sum.
Recall that the Steinberg module  $\St(G)$ of an arbitrary finite group $G$
is defined to be the virtual\break
$\Z_{(p)}G$-module
$$\St(G) = L(\St_*(G)).$$
(This definition of course coincides with the pre-existing definition for a
finite group of Lie type up to a sign.)
Note that this is a projective virtual\break $\Z_{(p)}G$-module, since its
reduction modulo $p$ is projective.

\proclaimtitle{cf.\ Webb {\cite[Thm.\ 6.8]{webb87arcata}}}
\proclaim{{C}orollary}
Let $G$ be a finite group and let $\cC$ be a collection of $p$\/{\rm
-}\/subgroups
 closed under passage to $p$-radical overgroups. Let $\fF: {\bO}(G)^{\op}
\to (\mbox{finitely generated } \Zpmod)$ be
any subgroup $\cC$\/{\rm -}\/acyclic functor{\rm .} Then
$$\fF(G) \cong - \sum_{P \in (\cC \cap \pRe)/G} \Hom_{WP}(\St(WP),\fF(P)).$$
If furthermore $CP$ acts trivially on $\fF(P)$ for all $P \in \cC$ then
$$\fF(G) \cong - \sum_{P \in (\cC \cap \pD)/G}
\Hom_{WP}(\St(NP/PC(P)),\fF(P)).$$
\endproclaim 

\demo{Proof} 
Let $\cF$ be the $G$-local coefficient system on $\cC$ given by precomposing $\frak F$
with the composite $\cC_G \to \bO_\cC \hookrightarrow {\bO}(G)$.
By Theorem~\ref{thm1} and our assumption we have that $\fF(G) =
H^0_G(|\cC|;\cF)$, and the higher cohomology is trivial. Hence $\fF(G) = L(
C^*_G(|\cC|;\cF))$.
Let $\fF_P$ denote the atomic functor associated to $\fF$, which is $F(Q)$
on subgroups $Q$ conjugate to $P$ and zero otherwise, and let $\cF_P$ denote
the corresponding $G$-local coefficient system. By filtering $\cF$ by the
functors $\cF_P$, as in the proof of Corollary \ref{cor1'} and using
the fact that this filtration of $C^*_G(|\cC|;\cF)$ splits as
$\Z_{(p)}$-modules, we get that
$$L(C^*_G(|\cC|;\cF)) \cong \sum_{[P] \in \cC/G} L(C^*_G(|\cC|;\cF_P)).$$
But the proof of Theorem~\ref{thm1} shows that $C^*_G(|\cC|;\cF_P)$ is
$\Z_{(p)}$-chain homotopy equivalent to $\Hom_{WP}(\Sigma\St_*(WQ),\fF(P))$
where $\Sigma$ here denotes
a  suspension of chain complexes. Combining these statements and using the
definition of $\St(G)$ we prove  the
corollary in the first case. The case where $CP$ acts trivially on $\fF(P)$
follows from the first statement by
Corollary~\ref{stcritcor} and Proposition~\ref{coinvprop}.
\enddemo 

\numbereddemo{{R}emark}
The above techniques can also be used to obtain information even if
acyclicity does not hold, although in that case the results often come as an
equality between virtual $\F_p$-modules, rather than as exact sequences. We
illustrate this with an example.

Suppose that $G$ is a finite group and $H$ is a normal subgroup. Then $G$
acts on ${\cal S}_p(G/H)$ by conjugation. Consider the Bredon equivariant
cohomology of  $|{\cal S}_p(G/H)|$ with values in the generic coefficient
system $\fF(-) = H^n(-;M)$, where $M$ is any $\F_pG$-module. By a spectral
sequence argument, as in the proof of Theorem~\ref{normsharp}, one sees that
$$\tilde H^*_G(|{\cal S}_p(G/H)|; H^n(-;M)) = \tilde H^*_G(|{\cal
S}_p(G/H)|; H^0(-;H^n(H;M))), $$
which need not equal zero. However, writing out the defining chain complexes
on the left- and the right-hand side
and rearranging we obtain
\begin{eqnarray*}
H^n(G;M) &\cong& \left(\sum_{[\sigma] \in |{\cal
S}_p(G/H)|/G}(-1)^{|\sigma|}H^n(G_\sigma;M)\right)  \\[6pt]
&&-\ 
\Hom_{G/H}(\St(G/H),H^n(H;M)) .
\end{eqnarray*}\break This result is due to Adem-Milgram \cite[Thm.\
2.2]{AM92}.
\enddemo

\numbereddemo{{R}emark} \label{detectionrem} Note that all our exact
sequences exist if just
$\lim^*_\D F=0$ for $*>0$, without making assumptions on $\lim^0_\D F$, if
one is willing to replace $\fF(G)$ with
$\lim^0_\D F$. This weaker assumption is satisfied more frequently in
nature. In fact, for {\it any} Mackey functor
$\fF: \bO(G)^{\op} \to \Zpmod$ we have  $\lim^*_{\bO_{\pSe}}F^{\beta} = 0$
for $*>0$ (see \cite{JM92}, \cite{JO96}).
The condition that $\fF(G) \to \lim^0_{\bO_{\pSe}}F^{\beta}$ is an
isomorphism is in practice more restrictive, and
requires $\fF$ to be projective relative to $p$-subgroups \cite{webb99}.
(Or, almost equivalently, satisfy an
`induction theorem' relative to $p$-subgroups \cite{webb99},
\cite[5.6]{benson91}.) For a Mackey functor $\fF$ which
does not satisfy this, one can instead consider the smallest collection
$\cC$ which is closed under passage to
subgroups, contains a collection relative to which $\fF$ is projective, and
which satisfies the condition of Remark
\ref{extendrem}. $\fF$ is now subgroup $\cC$-acyclic and Theorem~\ref{thm1}
still applies. We hope to examine some
interesting examples of this situation in detail in a future paper.
\enddemo
\pagebreak

\section{Sharpness} \label{sharpsec}

In this short section we explain how to obtain sharpness results from the
acyclicity results of the previous section. We will be brief and refer the
reader to \cite{dwyer98} for a much fuller treatment.

Let $M$ be some $\Z_{(p)}G$-module. A collection $\cC$ is called $M$-ample
if\break $H^*(G;M) \pil{\simeq}
H^*_{hG}(|\cC|;M)$. A collection $\cC$ is called subgroup $M$-sharp if
$H^n(-;M)$ is subgroup $\cC$-acyclic for all $n
\geq 0$. We define centralizer $M$-sharp and normalizer $M$-sharp similarly.
(We abbreviate $\Z_{(p)}$-sharp to just
sharp, and similarly for ample.)

Recall that for any $G$-space $X$ we have an isotropy spectral sequence
$$ E_2^{i,j} = H^i_G(X;H^j(-;M)) \Rightarrow H^{i+j}_{hG}(X;M).$$
The rewritings of Proposition \ref{prop1}, \ref{prop2} and \ref{normprop}
identify the $E_2$-term of the isotropy spectral sequence for $|E\bO_\cC|$,
$|E\A_\cC|$, and $|\cC|$ with the higher limits over $\bO_\cC^{\op}$,
$\A_\cC$, and $((\sd
\cC)/G)^{\rm op}$ respectively. Therefore the definition of\break
$M$-sharpness says that the $E_2$-term of the
isotropy spectral sequence must be concentrated on the vertical axis, where
we must have $E^{0,*} = H^*(G;M)$. In
particular, if $\cC$ is either subgroup, centralizer or normalizer $M$-sharp
then $\cC$ is $M$-ample, since
$|E\bO_\cC|_{hG} \simeq |\cC|_{hG} \simeq |E\A_\cC|_{hG}$. (Compare the
isotropy spectral sequence of $|E\bO_\cC|$,
$|E\A_\cC|$, or $|\cC|$ to that of a point.)

We can now use the acyclicity results in the previous sections to obtain
sharpness statements. We start by
demonstrating that $\pS$ is normalizer $M$-sharp if and only if  it is
$M$-ample. We have seen that the isotropy
spectral sequence for $|\pS|$ is concentrated on the two axes. On the
horizontal axis we have $E_2^{*,0} =
H^*_G(|\pS|;M)$, whereas the vertical axis is $E_2^{0,*} = H^*(G;M)$, for
$*>0$. Furthermore, this spectral sequence
has to collapse at
$E_2$. To see this note that  the map $\pS \to *$ has a $P$-equivariant
section given by $* \mapsto ZP$, and so   there
can be no differentials   in the isotropy spectral sequence for the
$P$-space $|\pS|$. But by a transfer argument, the
$E_2$-term for the $G$-space $|\pS|$ is a retract of the $E_2$-term for
$|\pS|$ viewed as a $P$-space; hence  there
can be no differential for the $G$-space either. Example~8.9 now gives us
the following result.

\proclaim{Theorem} \label{sharpnessresult}
 The collection $\pDeM$ of Example~{\rm 8.9} is subgroup\break $M$\/{\rm
-}\/sharp and normalizer $M$\/{\rm
-}\/sharp
 {\rm (}\/and is in particular
$M$\/{\rm -}\/ample{\rm ).}
 The collection $\pDM$ of Example~{\rm 8.9} is subgroup $M$-sharp if and
only if it is normalizer $M$\/{\rm
-}\/sharp if and only if it is
$M$\/{\rm -}\/ample if and only if $\Hom_G(\St_*(G),M)$ is an acyclic chain
complex{\rm .} 
\endproclaim 

Note that in the extreme case where $p \hbox{$\not |$}\   |G|$ we get that
$\pDM$ is\break $M$-ample if and only if
$M^G = 0$. The above result generalizes the result of Dwyer \cite{dwyer97},
\cite{dwyer98} that the collection of
$p$-centric and
$p$-radical subgroups is subgroup sharp, and the result originally due to
Webb \cite{webb91} that the collection of
$p$-radical subgroups is normalizer sharp.  The consequence that the
collection of $p$-radical and $p$-centric
subgroups is normalizer sharp was suspected by Smith-Yoshiara \cite{SY97}.
We remark that the analysis shows that
in the cases we have studied the collections lying in between our sharp
collections and $\pSe$ are also sharp.

\section{Stable elements and Alperin's fusion theorem}\label{stableelements}

In this section we explain how our work relates to Cartan-Eilenberg's
formula for stable elements in group
cohomology and Alperin's fusion theorem.

Let $\cC$ be a collection of $p$-subgroups of a finite group $G$, closed 
under passage to $p$-radical overgroups, and let
$F: \bO_\cC^{\op} \to \Zpmod$ be any functor. We call $\lim^0_{\bO_\cC}F$
the stable elements of $F$ with respect to $\cC$. (These are the stable
elements of Cartan-Eilenberg \cite{CE56} if $\cC = \pSe$ and $F =
H^n(-;M)$.) 

More amenable formulas for the stable elements can be obtained from
Alperin's fusion theorem. This is mentioned e.g.\  in Holt \cite{holt77},
and a related approach to the fusion theorem and stable elements is
presented in the pioneering work of Puig \cite[Ch.\ III \S 2]{puig76}.

Our point is that   $\lim_{\bO_\cC}^0F$ in our model for $\lim_{\bO_\cC}^*$F
gives the same kind of reduction in the
formula for stable elements, and the method gives some additional insight as
well. We start by explaining the results,
and will afterward relate them to Alperin's fusion theorem.

Let $\pH$ denote the collection of $p$-subgroups $Q$ in $G$ such that\break
${\cal S}_p(NQ/Q)$ is either empty or
disconnected (see Remark \ref{stronglypembedded}). In particular ${\cal
S}_p(NQ/Q)$ will be noncontractible so that
$\pH \subseteq \pRe$.  If $G$ is a finite group of Lie type of
characteristic $p$  then one easily identifies $\pH$
inside $\pRe$ as exactly the unipotent radicals of parabolic subgroups of
rank less than or equal to one.

\numbereddemo{{R}emark}
The role of the subgroups in $\pH$ for the $p$-local structure of $G$ has
been examined  especially by Goldschmidt
\cite{goldschmidt70} and Puig \cite{puig76}, and a link to subgroup
complexes was noted by Quillen \cite[\S
6]{quillen78}. (Puig calls the non-Sylow $p$-subgroups in $\pH$
$1$-essential (or just essential) \cite[p.\ 25]{puig76}
and non-Sylow subgroups in $\pH \cap \pD$ $C$-essential \cite[p.\
38]{puig76}; see also \cite{thevenaz95} where the
notation is modified.)
\enddemo

\numbereddemo{{R}emark} Easy examples show that $\pH$ is in general not an
ample collection.
 For example take $G = \GL_4(\F_2)$. If ${\cal H}_2(G)$ was an ample
collection it would be sharp and Corollary
\ref{exactliecor} should give an exact sequence e.g.\  for $F =
H^1(-;\F_2)$. However this is not the case, as can
perhaps most easily be seen by comparing it to the corresponding sequence
for the collection $\pR$.

Intuitively, one should think of the $G$-poset
 $\pH$ as giving generators and relations for the cohomology of $G$, whereas
the $G$-poset $\pR$ also encodes
relations amongst the relations necessary to get exact (sequence)
expressions for the cohomology of $G$.
\enddemo

Examination of  what happens at the $(0,0)$-spot of the spectral sequence of
Theorem \ref{sscor} gives the following
result which says that one can determine the stable elements `locally', 
i.e., via information coming from
$p$-subgroups and their normalizers; one does not have to consider
conjugation between different subgroups as it
appears from the definition.

\proclaim{Theorem} \label{absgenCE} Let $\cC$ be a collection of $p$\/{\rm
-}\/subgroups
in $G$ closed under passage to $p$\/{\rm -}\/overgroups{\rm ,}
 and let $F: \bO_\cC^{\op} \to \Zpmod$ be any functor{\rm .} Let $P$ be a
Sylow $p$\/{\rm -}\/subgroup in $G${\rm .} Then
$$ \lim_{\bO_\cC}{}^0F = \{ x \in F(P)^{WP} | \res^P_Q(x) \in F(Q)^{WQ}
\mbox{ for all } [Q] \in\cC \cap \pH/G\}$$
where $Q$ is an arbitrary representative of $[Q]$ such that $Q\leq P${\rm .}

If $CQ$ acts trivially on $F(Q)$ for all subgroups $Q \leq P${\rm ,}
 then   $\cC \cap \pH$ may be replaced by $\cC \cap \pH \cap \pD${\rm .}
\endproclaim 

\demo{Proof} 
Choose a height function on $\cC \cap \pH$ such that $\hi(P)=0$. Consider
the spectral sequence of Theorem \ref{sscor}. We have that $E_1^{i,-i} \neq
0$ only if $i=0$, since $\St_*(H)$ is $(-1)$-connected if $p | |H|$. Hence,
we can calculate $\lim_{\bO_\cC}^0F$ as the intersection of the kernels of
the differentials originating at the $(0,0)$-spot. By Proposition
\ref{h0formula} we have
\begin{eqnarray*}
E_1^{i,-(i-1)}& = &\bigoplus_{[Q] \in (\cC \cap \pRe)/G; \hi(Q) = -i}
\Hom_{WQ}(\tilde H_0({|{\cal
S}_p(WQ)}|),F(Q))\\[4pt]
& =& F(Q)^{\widetilde{WQ}}/F(Q)^{{WQ}}.
\end{eqnarray*}  In particular, only differentials going to a subgroup $Q
\in \cC
\cap \pH$ can contribute anything.

Finally, one easily verifies from the construction of the spectral sequence
(see the proof of Theorem \ref{sscor}) that
the kernel of the differential $d_n: E_n^{0,0} \to E_n^{n,-(n-1)}$ is just
equal to the kernel of $F(P)^{NP} \pil{\res^P_Q}
F(Q)
\to F(Q)/F(Q)^{{WQ}}$. This proves the first part of the theorem.

The last part follows since if $CQ$ acts trivially on $F(Q)$ then Theorem
\ref{sscor} holds with $\cC \cap \pRe$ replaced by $\cC \cap \pD$.
\enddemo

Specializing to the case of group cohomology, we find that the above result
takes the following form, in view of the results of
Section~\ref{exactseqsec}.
\pagebreak

\proclaim{{C}orollary}
\label{genCE} Let $M$ be any $\Z_{(p)}G$\/{\rm -}\/module{\rm ,} and let $P$
be a Sylow $p$\/{\rm -}\/subgroup in
 $G${\rm .} Then the following
formula for the group cohomology of $G$ with coefficients in $M$ holds\/{\rm
:}
\begin{eqnarray*}
H^*(G;M)& = &\{ x \in H^*(P;M)^{WP} | \res^P_Q(x) \in H^*(Q;M)^{WQ}\\
&&\quad \mbox{ for all } [Q] \in
(\pH\cap\pDeM)/G\},
\end{eqnarray*} where $Q$ is an arbitrary representative of $[Q]$ such that
$Q\leq P${\rm .}  
\endproclaim 

Notice also that Theorem \ref{olivercomplex} yields a dual version of
Theorem \ref{absgenCE}.

\proclaim{Theorem}\label{dualgenCE} Let $\cC$ be any collection of $p$\/{\rm
-}\/subgroups 
of $G$ closed under passage to
nontrivial elementary abelian $p$\/{\rm -}\/subgroups{\rm ,} and let $F :
\A_{\cC} \to \Zpmod$ be any functor{\rm .}
 Then
\begin{eqnarray*}\lim_{\A_{\cC}}{}^0F& = &\Big\{ \prod x_V \in
\hspace{-20pt}\prod_{[V] \in \cC'/G, \rk V =1}
\hspace{-20pt} F(V)^{NV} | \sum_{[V],[V] \leq [W]} \iota(x_V) \in F(W)^{NW}
\\
&&
\hbox{ for all } [W] \in \cC'/G, \rk W =2
\Big\}.
\end{eqnarray*}    Here $\cC' = \cC \cap \pA$ and $\iota$ is the inclusion
map $V \leq W${\rm ,}
 where   representatives $V$
such that $V \leq W$ are chosen{\rm .}
\endproclaim  

For group cohomology Theorem \ref{dualgenCE} gives an expression for
$H^*(G;M)$ in terms of centralizers of
elementary abelian $p$-subgroups of rank $1$ and $2$.

\vglue12pt {\it Example} 10.6  (The extra-special group $p^{1+2}_+$).
Consider a finite group $G$
with Sylow
$p$-subgroup the extra-special group $P = p^{1+2}_+ = \langle a,b,c |
a^p=b^p=c^p=1, [a,c]=[b,c]=1,[a,b]=1 \rangle$. If $Q
\leq P$ satisfies $C_PQ = ZQ$ then $Q=P$ or $Q \cong \Z/p \times \Z/p$.
Furthermore, if $Q \cong \Z/p \times \Z/p$
and $W = NQ/CQ \leq \Out(Q) = \GL_2(\F_p)$ satisfies $O_p(W) = e$ then an
easy computation reveals   that
$\SL_2(\F_p) \leq W$, since $p | |W|$. Therefore
$$\pD/G = \{[P],[V_1], \ldots , [V_k]\}$$ where $V_1, \ldots, V_k$ are
$G$-conjugacy class representatives of the rank $2$ elementary abelian
subgroups of $P$ with the property that $\SL_2(\F_p) \leq NV_i/CV_i$.

Corollary \ref{genCE}, when applied to this example, gives a generalization
of \cite[Thm.\ 0.2]{YT96}, where a formula of this type was proved for
certain sporadic groups. The fact that the collection $\pD$ is ample
(Theorem \ref{sharpnessresult}) likewise generalizes \cite[Thm.\
1.8]{yagita99}.
\vglue12pt

{\it Example} 10.7 (Swan groups).   Suppose that for a $p$-group $P$ the
condition
 $$C_P(Q) = ZQ \mbox{ and } N_P(Q)/Q \cap O_p(\Out(Q)) = e,$$
mentioned \pagebreak in Remark \ref{localconditions},
is not satisfied for any proper subgroups $Q$ of $P$. Then we conclude that
for {\it any} finite group $G$ having $P$ as Sylow $p$-subgroup $\pD/G =
[P]$ and hence $H^*(G;\F_p) \cong H^*(P;\F_p)^{WP}$,
by Corollary \ref{genCE}.
 This was originally proved by Martino-Priddy \cite[Thm.\ 2.3]{MP97}.
\advance\theoremcount by 2

\numbereddemo{{R}emark} \label{mislinrem} Using Goldschmidt's or Puig's
version of Alperin's fusion theorem (see below) one sees that $\pD/G =[P]$
for the Sylow $p$-subgroup $P$ of $G$
if and only if  $NP$ controls (strong) $p$-fusion in $G$.
Furthermore,\break Mislin's theorem \cite{mislin90} tells us
that 
$H^*(G;\F_p)
\pil{\simeq} H^*(P;\F_p)^{WP}$ implies that $NP$ controls (strong)
$p$-fusion in $G$, so that the three conditions
$\pD/G =[P]$,   $NP$ controls (strong) $p$-fusion in $G$,  and $H^*(G;\F_p)
\pil{\simeq} H^*(P;\F_p)^{WP}$, and are in
fact equivalent. In general, one can interpret the number $\dim |\pD|$ as a
measure of how complicated the $p$-fusion
is in
$G$.
\enddemo

\numbereddemo{{R}emark} The proof of Theorem~\ref{absgenCE} in the case
where $CQ$ acts trivially on $F(Q)$ allows for the seemingly stronger
statement that we may restrict attention to $p$-centric subgroups $Q \in
\cC$ for which ${\cal S}_p(NQ/QC(Q))$ is empty or disconnected. However,
using some group theory one can show that this is equivalent to assuming $Q
\in \cC \cap \pH \cap \pD$
(see \cite[Ch.\ III,  Cor.~2]{puig76}).
\enddemo 

\numbereddemo{{R}emark} Note that the proof of Theorem~\ref{absgenCE}
actually shows more than stated in the
theorem. We already know that  $\res^P_Q(x) \in F(Q)^{\widetilde{WQ}}$, so we
only have to calculate invariants under a
subgroup which together with $\widetilde{WQ}$ generate $WQ$. Also, the proof
shows that each $Q \in \cC \cap \pH$ can
only contribute if $\Hom_{WQ}(\tilde H_0({{\cal S}_p(WQ)}),F(Q)) \neq 0$.
Remark \ref{lim1formula} gives us a good
handle on when this happens, and can allow us to rule out more subgroups if
we know something about $F(Q)$. This can
for instance be useful for low dimensional group cohomology calculations, as
we shall see in the next example.
\enddemo

\numbereddemo{Example} \label{miyamoto}Suppose that $G$ is a finite group,
whose Sylow $p$-subgroup $P$ has nilpotency class at most $c$, i.e., any
$c$-fold commutator between elements in $P$ is trivial.
If $Q \leq P$ and $g \in WQ$ of order $p$, then Miyamoto has proved
\cite[Lemma 2]{miyamoto81}, using induction and the Lyndon-Hochschild-Serre
spectral sequence, that
$(1-g)^{(c-1)n+1}$ acts as zero on the $WQ$-composition factors of
$H^n(Q;\F_{p})$. Now, if $(c-1)n+1 \leq p-1$ then, by Example
\ref{st0refinement}, $\Hom_{WQ}(\tilde H_0({{\cal S}_p(WQ)}),H^n(Q;\F_{p}))
= 0$ for all $Q < P$. Hence, $$\mbox{if } c \leq \frac{p-2}n + 1 \mbox{
then } H^n(G;\F_{p}) = H^n(P;\F_{p})^{NP}.$$
This result, due to Miyamoto \cite{miyamoto81}, generalizes both Swan's
classical theorem, which is the case $c=1$, as well as results for $n=1$ and
$2$ of Glauberman and Holt \cite{holt77} respectively.
\enddemo

\numbereddemo{{R}emark} \label{fusionrelaterem} The conditions on subgroups
are stated in a number of different ways in the literature. For the
convenience of the reader, we state some relations between these.
\begin{itemize}
\item[1.] \label{Rone} A subgroup $Q$ of $P$ is a {\it tame intersection} in
$P$ if there exists a Sylow $p$-subgroup
$P'$ of $G$ such that $Q=P \cap P'$ and that $N_P(Q)$ and $N_{P'}(Q)$ are
both Sylow $p$-subgroups in $NQ$. Puig
observed \cite[p.\ 42]{puig76} (see also \cite[\S 6]{quillen78}) that if $Q
\in \pH$ and if $N_P(Q)$ is a Sylow
$p$-subgroup in $NQ$, then such $P'$ always exists. Hence, by Sylow's
theorem, any $Q \leq P$ with $Q \in \pH$ is
$G$-conjugate to a subgroup $Q'\leq P$ with $Q'$ a tame intersection in $P$.

\item[2.] A subgroup $Q \leq P$ is $p$-centric in $G$ if and only if  for
all subgroups $Q' \leq P$ with $Q'$
$G$-conjugate to
$Q$ we have $C_P(Q') = ZQ'$. (Use Sylow's theorem.)

\item[3.] If $Q$ is assumed to be a tame intersection in $P$ and $C_P(Q) =
ZQ$ then $Q$ is $p$-centric. To see this just
note that since $N_P(Q)$ is assumed to be a Sylow $p$-subgroup in $NQ$ and
$CQ \unlhd NQ$, $C_P(Q) = CQ \cap
N_P(Q)$ is also a Sylow $p$-subgroup in $CQ$.

\item[4.] Note that we always have $O_{p'}(NQ) = O_{p'}(CQ)$.
Hence, if $Q$ is $p$-centric then the condition that $O_p(NQ/QC(Q)) = e$ is
equivalent assuming that $Q$ is a Sylow $p$-subgroup of $O_{p',p}(NQ)$. If
$Q$ is assumed to be a tame intersection in $P$, it is again equivalent to
assuming $O_{p',p}(NQ) \cap P = Q$.
\end{itemize}
By 2 we can view the $p$-centric condition as being a $G$-conjugacy
invariant refinement of the condition $C_P(Q) =
ZQ$. Also,  1 and 3 show that the tame intersection condition picks out
distinguished representatives from a
$G$-conjugacy class of $p$-centric subgroups. (The name $p$-centric subgroup
was coined by Dwyer \cite{dwyer97},
inspired the $p$-centric diagrams of \cite{DK92}. In parts of group theory
$p$-centric subgroups are called
self-centralizing \cite[p.\ 324]{thevenaz95}.)
\enddemo

\numbereddemo{Definition} \cite{alperin67} Let $G$ be a finite group and $P$
a fixed Sylow\break $p$-subgroup in
$G$. A set ${\Bbb X}$ of pairs $\{(Q,T)\}$, where $Q$ is a subgroup of $P$
and $T$ is a subset of $NQ$, is called a {\it
conjugation family} provided the following condition is satisfied: For any
subset $A$ of $P$ and any $g \in G$ such
that $A^g \leq P$, there exist pairs $(Q_1,T_1),\ldots, (Q_n,T_n)$ in
${\Bbb X}$ and elements $x_1,\ldots,x_n,y$ of $G$ such that $g = x_1\cdots
x_ny$, $x_i \in T_i$ for $i=1,\ldots,n$, $y
\in NP$, and $A^{(x_1\cdots x_i)} \leq Q_{i+1}$ for $i=0,\ldots,n-1$.  ($A^g
=g^{-1}Ag$.)
\enddemo
Note that  for each conjugation family ${\Bbb X}$ we  can form a new
conjugation family replacing each $(Q,T)$ with
$(Q,\bar T)$, where $\bar T$ is the subgroup generated by the union of the
subsets $\{ U |  (Q,U) \in {\Bbb X}\}$, which
will be just as good for practical purposes.  Hence we can when convenient
assume that $(Q,T)$ is uniquely determined
by $Q$, and that $T$ is a group.
 We say that a conjugation family is {\it invariant} if it, in addition to
the just stated requirements, satisfies the fact
that if
$(Q,T) \in {\Bbb X}$ and $Q^g \subseteq P$ then $(Q^g,T^g) \in {\Bbb X}$.

It is  not obvious that conjugation families do indeed exist. However the
following is immediate.

\proclaim{Proposition} \label{conjfamstable}
 Let $\cC = \pSe$ and let $F: \bO_\cC^{\op} \to \Zpmod$ be any functor{\rm
.} Assume that ${\Bbb X}$ is a conjugation
family in $G${\rm .} Then
$$\lim_{\bO_\cC}{}^0F = \{ \alpha \in F(P)^{WP} | \res^P_Q(\alpha) \in
F(Q)^{T}  \mbox{ for all } (Q,T) \in {\Bbb X}\} .$$

If ${\Bbb X}$ is invariant then it is enough to check the restriction
condition on one representative $(Q,T)$
for each $G${\rm -}conjugacy class of pairs in ${\Bbb X}${\rm .}
\endproclaim 
\demo{Proof} 
The first part of the statement follows by writing out the definition of
$\lim^0$ (see e.g.\  \cite{holt77}). The last part of the statement follows
by a downward induction on the size of $Q$. (The statement is trivial for
$Q=P$.)
\enddemo 

Goldschmidt's version
 of Alperin's fusion theorem \cite{alperin67}, \cite{goldschmidt70} (stated
in a stronger form due to Miyamoto
\cite[Cor.\ 1]{miyamoto77}, and slightly reformulated) says that the set of
pairs $\{(Q,T)\}$ where $Q$ is a tame
intersection subgroup in $P$ such that $Q \in \pH$ and $T = NQ$ if $Q \in
\pD$ and $T = CQ$ otherwise, is a
conjugation family. Hence, Proposition \ref{conjfamstable} combined with the
fusion theorem produces a statement
very close to Theorem~\ref{absgenCE}. (Compare also \cite[Ch.\ III \S
2]{puig76}.)


\bye